\documentclass[12pt]{amsart}\usepackage{a4wide,enumerate,color,graphicx}
\usepackage{amsmath}
\usepackage{scrextend} 
\usepackage{amssymb}
\usepackage{url,hyperref}
\usepackage{amsfonts}
\usepackage{amsthm}
\usepackage{enumerate}
\usepackage{mathrsfs}
\usepackage{caption}
\usepackage{mhchem}     

\newtheorem{lemma}{Lemma}[section]
\newtheorem{prop}[lemma]{Proposition}
\newtheorem{theo}[lemma]{Theorem}

\newtheorem{rem}[lemma]{Remark}
\newtheorem{coro}[lemma]{Corollary}

\DeclareMathOperator{\divv }{div}

\DeclareMathOperator{\diff }{\delta}

\DeclareMathAlphabet{\mathdutchcal}{U}{dutchcal}{m}{n}

\allowdisplaybreaks

\begin{document}

\title[Transmission problem near a cusp]{Gradient and Hessian regularity in elliptic transmission problems near a point cusp}


\author[D.~Bothe, P.-E.~Druet, R.~Haller]{Dieter Bothe, Pierre-Etienne Druet and Robert Haller}
\address{Institute for Mathematics, TU Darmstadt}
\email{bothe@mma.tu-darmstadt.de,druet@mma.tu-darmstadt.de,\newline haller@mathematik.tu-darmstadt.de}

\date{\today}

\begin{abstract}
We consider elliptic transmission problems in several space dimensions near an interface which is $C^{1,1}$ diffeomorphic to an axisymmetric reference-interface with a singular point of cusp type. 
We establish the regularity of the gradient and of the Hessian in $L^p$ spaces up to the cusp point for local weak solutions. We obtain regularity thresholds which are different according to whether the cusp is inward or outward to the subdomain, and which depend explicitly on the opening of the interface at the cusp. Our results allow for source terms in the bulk and on the interface.
\end{abstract}
\maketitle

\setcounter{tocdepth}{2}
\tableofcontents

\section{Introduction}

We consider the regularity of solutions to the elliptic interface problem
\begin{align}
\label{P1}	-\divv(\kappa(x) \, \nabla u) = f(x) \quad \text{ in }\quad \Omega \setminus S\, ,\\
\label{P2}	[\![-\kappa(x) \, \nabla u]\!] \cdot \nu(x) = Q(x) \quad [\![u]\!] = 0 \quad \text{ on } S \, ,
	\end{align}
posed in a domain $\Omega \subset \mathbb{R}^d$, $d \geq 2$, with embedded hyper-surface $S \subset \Omega$. We assume that $S$ is regular up to a single point $\{x^0\}$ in which a singularity of cusp-type occurs. In \eqref{P2} the double bracket stand for the jump of a quantity across the surface $S$. Our aim is to prove the existence and integrability properties for second derivatives \emph{up to the singular point} for local weak solutions $u$ to \eqref{P1}, \eqref{P2}. Moreover, we shall discuss the regularity of $\nabla u$ in $L^p$ spaces for $p > 2$. Notice that, in the present context, {\it local regularity} is always meant as a property valid in a neighbourhood $B_{R_0}(x^0)$ and up to $x^0$. Consequently, it is not necessary to further specify the regularity of the domain $\Omega$, or to formulate far-field conditions for \eqref{P1}, \eqref{P2}. 

In this paper, we restrict to discussing an axial symmetric reference configuration. Fixing for simplicity $x^0 = 0$, we hence assume that there is a non-negative, non-decreasing function $\sigma\in C^2(]0,R_0]) \cap C([0,R_0])$ for a $R_0 > 0$ such that $\sigma(0) = 0$ and 
\begin{align}\label{cusp1}
	\liminf_{\rho \rightarrow 0+} \sigma^{\prime}(\rho) \in ]0,\,  +\infty] \, ,
\end{align}
such that, possibly after an affine transformation of the coordinate system, the surface $S$ possesses the graph-representation
\begin{align}\label{itisagraph}
S \cap \Omega := \{x = (\bar{x}, \, x_d)\in \Omega \, : \, x_d = \sigma(|\bar{x}|)\} \, .
\end{align}
Here $\bar{x} = (x_1, \ldots, x_{d-1})$ denotes generic elements of $\mathbb{R}^{d-1}$.

Under these conditions, $S$ is oriented and it partitions $\Omega$ into two domains, which we call $\Omega_i$, $i = 1,2$. We here adopt the convention that $\nu = \nu(x)$ occurring in \eqref{P2} is the unit normal vector to $S$ pointing into $\Omega_2$ located above the surface $S$, which we also call the \emph{outward} cusp domain. Since $\sigma$ is of class $C^2$, and since we are interested in local regularity near $x^0 = 0$, there is in view of \eqref{cusp1} no loss of generality to further assume that $\sigma$ is strictly increasing on $]0, \, R_0[$. Otherwise there is a smaller domain of definition $]0,R_1[$ where $\sigma^{\prime}$ is positive. 

\begin{minipage}{\textwidth}
\centering
	\includegraphics[width=0.6\textwidth]{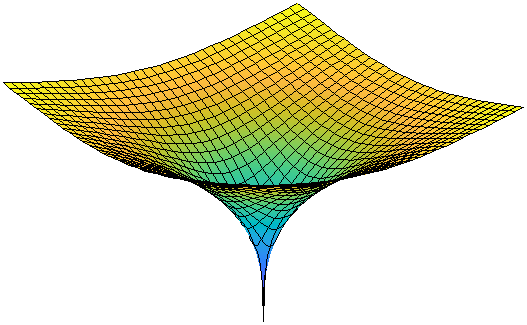}
	\captionof{figure}{The graph of $\sigma(r) = r^{\frac{1}{3}}$ in the unit square.}
\end{minipage}

 As we will show, our results are valid for $C^{1,1}-$transformations of the reference configuration, yielding more generality. For instance, the \emph{outward} cusp-domain $G$ considered in \cite{maz2010solvability}, \cite{nazarov2016} is characterised, in an open neighbourhood $U$ of the origin, as $G \cap U = \{x = (\bar{x}, \, x_d) \, : \, (1/\phi(x_d)) \, \bar{x}  \in \omega\}$, where $\bar{x} = (x_1, \ldots,x_{d-1})$, $\phi$ is a $C^2$ increasing function with $\phi(0) = 0$ and $\phi^{\prime}(0) = 0$ and $\omega \subset \mathbb{R}^{d-1}$ is a bounded domain of class $\mathcal{C}^2$ containing the origin. For $\omega = B^{d-1}_1(0)$ (the unit ball centred at the origin in dimensions $d-1$), this is obviously the axisymmetric case. In the general case, we can choose a $C^{2}-$diffeomorphism that maps $\omega$ onto $B^{d-1}_1(0)$ in order to come back to our reference configuration. \\

We turn to considering the remaining data occurring in \eqref{P1} and \eqref{P2}. For $x\in \Omega\setminus S$, the matrix $\kappa(x)$ is assumed real-valued, symmetric, uniformly elliptic and bounded with respect to $x$. We introduce
\begin{align}\label{ellipt}
	0 < k_0 := \inf_{x\in \Omega \setminus S} \inf_{\eta \in \mathbb{R}^d, \, |\eta|=1} \kappa(x) \eta \cdot \eta \quad \text{ and } \quad +\infty > k_1 := \sup_{x\in \Omega \setminus S} \sup_{\eta \in \mathbb{R}^d, \, |\eta|=1} \kappa(x) \eta \cdot \eta \, .
\end{align}
Moreover, we assume that $f$ in the equation \eqref{P1} and $Q$ in \eqref{P2} are given locally integrable functions. In \eqref{P2}$_2$, we require the continuity of $u$ across the interface, which is the natural condition in the context of heat transfer problems where $u$ is e.g.\ the temperature. In the context of mass transfer, $u$ represents a concentration and a discontinuity of the type $u_{\Omega_2} = H \, u_{\Omega_1}$ occurs on $S$ with constant coefficient $H >0$, see \cite{BOTHE2013283}. This case can also be handled in the context of \eqref{P1}, \eqref{P2} by rescaling $u$ and $\kappa$ on the side of $\Omega_1$. Likewise, a prescribed discontinuity of $[[u ]]_S = u_S$ across the interface can be reduced to the homogeneous case using additive modifications of $f$ and $Q$.

These considerations allow us to define the concept of a \emph{local weak solution} to \eqref{P1}, \eqref{P2} as a function $u \in W^{1,p}_{\text{loc}}(\Omega)$ for a $1\leq p \leq +\infty$ such that $\int_{\Omega} \kappa \, \nabla u \cdot \nabla \phi \, dx + \int_S Q \, \phi \, da = \int_{\Omega} f \, \phi \, dx$ for all $\phi \in \mathscr{D}(\Omega)$. The conditions under which such a definition meakes sense will become clear below. Here, we abbreviate $da = d\mathscr{H}^{d-1}$ with the $(d-1)$-dimensional Hausdorff--measure $\mathscr{H}^{d-1}$, which is a surface measure on $S$ (area).\\

{\bf State of the art, motivation and originality of our contribution:} Domains with cusps and formation of cusp-singularities can be observed in several areas of physics. In fluid mechanics, the study of cusp-shapes is relevant for flows near objects deposited on a substrate (\cite{nazarov1995-2}), but they also occur in many applications with free surface driven by surface tension: For instance in connection with drops spreading on substrates \cite{karpitschka2015droplets}, or bubbles rising in non-Newtonian liquids (\cite{liu1995two,zenit2018hydrodynamic,bothe2022molecular}). More precisely, the cusp-shaped bubble occurs during the quasi-stationary rise of a supercritical volume of gas in a polymer solution. 
In the latter cases, cusps occur as singularities in curvature and may precede a topological change, like the so-called \emph{tip streaming} phenomenon, see \cite{siegel2000cusp,krechetnikov2015cusps,montanero2020dripping}.  The papers \cite{jeong1992free,shikhmurzaev1998cusped,shikhmurzaev2005singularity} asks under which type of conditions cusp-singularities are compatible with surface tension dynamics.
In the case of a cuspidal bubble or droplet with mass or heat transfer across the fluid interface, problems of type \eqref{P1}, \eqref{P2} result in the quasi-steady case, see \cite{eggers2015singularities}.

Cusps also occur as special singular flow profiles in gas dynamics. The recent paper \cite{ShkollerVicol} interpret the occurrence of these profiles as an indicator for upcoming shock formation. In solid mechanics, cusps are observed in connection with cracks (\cite{nazarov2016}) and also in form-optimisation (\cite{LiuNeittaanmaekiTiba2003}). A nice recent investigation concerns the morphogenesis of an apple (\cite{chakrabarti2021cusp}).

 While cusp singularities in two space-dimensions necessarily occur as points, in 3D one has to distinguish between cusp-points, like for the case of the tip streaming or the apple cusp, and cusp-edges, like in the case of a liquid spreading on a surface with zero contact angle. The case of a cusp-edge occurs quite frequently in applications, whenever the boundaries of two objects meet tangentially. Surprisingly, it can also occur for rising bubbles in polymer solutions, see \cite{hassager1979negative,ohta2015dynamic}. 
 
 Accordingly, domains with cusp-points or -edges have always been a relevant challenging topic for the regularity theory of PDEs. In particular, domains exhibiting an outward cusp do not satisfy the cone condition, even in its weaker forms (see \cite{adamsfournier}), so that many propositions of Sobolev space theory might fail, resp.\ have to be re-checked. 
 
 Concerning domains with boundary cusps in connection to elliptic PDEs, there are several studies. Some among the most interesting contributions are \cite{nazarov1995} on the Neumann--problem for second-order elliptic operators near cusp-points in 3D, \cite{steux1997} and \cite{dauge1996} concerning the Dirichlet problem for elliptic operators of order $2m$ in 2D. For more general approaches, we can mention \cite{kozlov1997elliptic}, Part 3, or \cite{bschulze}. All these papers develop a general operator--adapted setting to study elliptic problems in domains with cusps by means of certain weighted $H^s-$spaces. There are also investigations concerning the Navier-Stokes equations in 2D cusp-domains, see \cite{nazarov1995-2} or more recently \cite{kaulakyte2021nonhomogeneous}.

\emph{Interface} or transmission problems at cuspidal interfaces for elliptic equations seem to have been comparatively less studied. In \cite{nicaise}, two- and three-dimensional interface problems are investigated near cusp-edges. Here, the interface is smooth while the cusp is located on the outer boundary. The authors establish some piecewise $H^2-$regularity for the solution. In \cite{laurencot}, which treats the 2D case, the cusp likewise occurs in the domain due to the tangential junction of an interface with the outer boundary. Here too, the authors prove the piecewise $H^2-$regularity of the solution. In both investigations, the cusp-singularity occurs at an outer boundary where homogeneous Dirichlet conditions are assumed -- either directly or after subtracting suitable data. In such a case, it is possible to gain the $H^2-$regularity via the ad-hoc assumption of compatible boundary data: On this topic we refer to \cite{druetboheme} for the general investigation of the junction of an interface with the outer boundary in 3D (see also \cite{lemrabet77} for particular results). Completely different is the case of an internal junction of three or more interfaces (triple contact line), where no $H^2$ regularity is to be expected, as it would impose the complete vanishing of the gradient in some weak sense (see \cite{druet2018regularity}, or the former \cite{kellogg1971singularities}, for results in standard, non-weighted Sobolev--classes). 

It seems that at the exception of the investigation \cite{maz2010solvability}, dealing with harmonic extensions for surface-sources in dual spaces by means of singular integral operators, elliptic transmission problems near an internal cusp-point were not studied in the literature. This is precisely the subject of the present paper. We are concerned with the regularity of second generalised derivatives and the higher integrability of the gradient up to the singular interface for distributional solutions to the problem \eqref{P1}, \eqref{P2}, in general space dimensions. 
Our method is thoroughly elementary: We prove the existence of second derivatives in a direction determined by the opening of the cusp (see Section \ref{direction}), then study the $(d-1)-$dimensional transmission problems orthogonal to this direction, which exhibit a regular interface (Section \ref{dmoinsun}). Our estimates depend on two parameters, associated with the diffusion coefficient $\kappa$ and the space dimension: The Meyers constant $p_0 > 2$, which is the threshold for which the operator $\divv (\kappa \nabla \cdot)$ is bounded invertible from $W^{1,p}$ into $(W^{1,p^{\prime}})^*$, and the deGiorgi-Nash exponent $\alpha_0$ of H\"older regularity for local solutions to $\divv (\kappa \nabla u) = 0$. Moreover, in the region $\Omega_2$ exhibiting the outward cusp, they depend on its opening measured by the largest constant $0 < \theta \leq 1$ such that 
\begin{align*}
	\sigma^{\prime\prime}(\rho) \geq (\theta - 1) \, \sigma^{\prime}(\rho) /\rho \quad \text{ for all } \quad 0 < \rho \leq R_0 \, .
\end{align*}
For smooth data $f$ and $Q$, we obtain in Theorems \ref{main} and \ref{main2} the integrability thresholds $\min\{p_0^*, \, d/(1-\alpha_0)\}$ for the gradient in the domain with inward cusp, and in the domain with the outward cusp the threshold $\max\{p_0, \, \min\{p_0^*, \, (d+\theta-1)/(1-\alpha_0)\}\}$. Here $p_0^*$ is the Sobolev embedding exponent for $p_0$. For the Hessian, the thresholds are $\min\{p_0, \, d/(2-\alpha_0)\}$ in $\Omega_1$, and $\min\{p_0, \, (d+\theta-1)/(2-\alpha_0-(\theta-1)/d)\}$ in $\Omega_2$. 
 
For $\Omega_2$, these results partly rely on a Gagliardo-Nirenberg inequality for domains with cusps (Sections \ref{gnsection}, \ref{sources}). \\

{\bf Notation and assumptions:} We at first formulate all conditions that will occur for the representation $\sigma$ of the surface $S$ (cf.\ \eqref{itisagraph}):
\begin{labeling}{(A1)}
	\item[(A1)] $\sigma\in C^2(]0,R_0]) \cap C([0,R_0])$ satisfies $\sigma(0) = 0$ and $0 < \sigma^{\prime}(\rho)$ for all $0 < \rho \leq R_0$;
\item[(A2)] $\liminf_{\rho \rightarrow 0+} \sigma^{\prime}(\rho) = +\infty$ for a cusp;	
\item[(A2')] $\liminf_{\rho \rightarrow 0+} \sigma^{\prime}(\rho) = \sigma^{\prime}_0 >0$ for a cone;	
\item[(A3)] There is $0< \theta \leq 1$ such that $\sigma^{\prime\prime}(\rho) \, \rho/\sigma^{\prime}(\rho)  \geq \theta -1 $ for all $0< \rho\leq R_0$;
\end{labeling}
The additional condition (A3) measures the opening of the cusp, and it shall be used to obtain uniform estimates. A standard example is $\sigma(\rho) = \rho^{\theta}$, which satisfies (A1), (A2) and (A3) if $0 < \theta < 1$ and (A1), (A2'), (A3) if $\theta = 1$. A function satisfying (A1) and (A2) but not (A3) is $\sigma(\rho) = 1/\ln (R_0/\rho)$ for $0 < \rho \leq R_0$.

For $1 \leq p \leq +\infty$, we make use of the standard Lebesque spaces $L^p(\Omega)$ and Sobolev spaces $W^{k,p}(\Omega)$, $k = 1,2$. Moreover, we make use of the spaces $C^{\alpha}(\Omega)$ for $\alpha \in ]0,1]$. The definition and the standard topologies on these spaces are assumed to be known.

Before stating our theorems, we introduce some abbreviations. For $1 \leq p \leq + \infty$, we define $p^{\prime} = p/(p-1)$ to be the conjugated exponent ($1^{\prime}= +\infty$ and $+\infty^{\prime} = 1$). We denote by $p^* = p^*(d)$ the Sobolev embedding-exponent:
\begin{align*}
	p^* := \begin{cases}
		dp/(d-p) & \text{ for } 1 \leq p < d \, ,\\
		\delta^{-1}< + \infty \text{ with arb. } 0 < \delta <1 & \text{ for } p = d \, ,\\
		+\infty & \text{ for } p > d \, . 
	\end{cases}
\end{align*}
For commodity we write $p^{\rm \#}$ instead of $((p^{\prime})^*)^{\prime}$, which is the index
\begin{align*}
	p^{\rm\#} = \begin{cases}
		1 & \text{ for } 1\leq p < d/(d-1) \, ,\\
		1+ \delta \text{ with arb. } \delta>0 & \text{ for } p = d/(d-1) \, ,\\
		d\, p/(d+p)& \text{ for } p > d/(d-1) \, .
	\end{cases}
\end{align*}
For the bulk-to-surface Sobolev embedding exponent ($p_S^* = (d-1)p/(d-p)$ for $d < p$, $p_S^* $ arbitrary finite if $d=p$ and $p_S^* = + \infty$ otherwise), the exponent $p_S^{\rm \#} := ((p^{\prime})^*_S)^{\prime}$ is
\begin{align*}
	p^{\rm \#}_S:= \begin{cases}
		1 & \text{ for } 1 \leq p < d/(d-1) \, ,\\
		1+\delta \text{ with arb. } \delta >0  & \text{ for } p= d/(d-1) \, ,\\ 
		(d-1) \, p/d & \text{ for } p > d/(d-1) \, .  
	\end{cases}
\end{align*}
In order to formulate our result, we also need a notion for the trace of Sobolev functions on $S$. For $1 \leq p \leq \infty$, every element $u$ of the Sobolev--space $W^{1,p}(\Omega)$ possesses a trace on $S \setminus B_{\epsilon}(0)$ for all $\epsilon >0$. Let us call ${\rm trace}_{S \setminus B_{\epsilon}(0)}$ the corresponding trace operator mapping from $W^{1,p}(\Omega)$ to $W^{1-1/p}_p(S \setminus B_{\epsilon}(0))$. We then can define a function
\begin{align}
	v(x) := \Big({\rm trace}_{S \setminus B_{|x|/2}(0)}(u)\Big)(x) \quad \text{ for } \quad x \in S \setminus \{0\} \, ,
\end{align}
and we call it the trace of $u$ on $S$, with notation $v = {\rm trace}_S(u)$. Reversely, for a function $v \in L^1(S)$ and for $1 \leq p \leq \infty$, we shall write $v \in {\rm tr}_S(W^{1,p}(\Omega))$ if there is $u \in W^{1,p}(\Omega)$ such that ${\rm trace}_S(u) = v $. The norm is given by
\begin{align*}
	\|v\|_{{\rm tr}_S(W^{1,p}(\Omega))} := \inf_{u \in W^{1,p}(\Omega) \, : \, {\rm trace}_Su = v} \|u\|_{W^{1,p}(\Omega)} \, .
\end{align*}
We note that $S$ is a $d-1$ set in the sense of Jonsson and Wallin (see \cite{MR0820626}, pages 28--30), where we can choose the surface measure as the standard Hausdorff-measure $\mathscr{H}^{d-1}$. For $1< p < +\infty $ it is thus possible to characterise ${\rm tr}_S(W^{1,p}(\Omega)) = B^{p,p}_{1-\frac{1}{p}}(S)$ (see  \cite{MR0820626}, VII, Theorem 1, page 182), where the latter Besov space is characterised via
	\begin{align}\label{besovraum}
B^{p,p}_{\alpha}(S) = \Big\{v \in L^p(S; \, da) \, : \, \int\int_{|x-y|<1} \frac{|v(x) - v(y)|^p}{|x-y|^{d-1+\alpha \, p}}
da(x) da(y) < + \infty\Big\} \, ,		
	\end{align}
	see \cite{MR0820626}, page 103, (V.1.1). We refer also to \cite{MR3573649}, Proposition 2.5.
	
	In order to formulate assumptions for $\kappa$ and $f$, we further need to recall two basic results of elliptic theory that are crucial in the present context. We denote by $A$ the weak elliptic operator $-\divv (\kappa \, \nabla \cdot )$, then
	\begin{enumerate}[(I)]
		\item\label{Meyers} There is $p_0 > 2$ depending only on $d$ and the numbers $k_0$ and $k_1$ of \eqref{ellipt} such that for all $ p \in ]p^{\prime}_0, \, p_0[$ with $p^{\prime}_0 := p_0/(p_0-1)$, the weak elliptic operator $A : \, W^{1,p}_0(B_{R_0}(x^0)) \rightarrow [W^{1,p^{\prime}}(B_{R_0}(x^0))]^*$ is invertible with $A^{-1} \in \mathscr{L}([W^{1,p^{\prime}}(B_{R_0}(x^0))]^*, \, W^{1,p}_0(B_{R_0}(x^0)))$;
		
		\item\label{deGiorgiNash} There is $0 < \alpha_0 \leq 1$ depending only on $d$ and the numbers $k_0$ and $k_1$ of \eqref{ellipt} such that for $p >d$, $A^{-1}$ is linear and bounded from $[W^{1,p^{\prime}}(B_{R_0}(x^0))]^*$ into the H\"older space $C^{\min\{\alpha_0,1-d/p\}}(B_{R_0}(x^0))$.
	\end{enumerate}
	The exponent of \eqref{Meyers} was first obtained in \cite{meyers} (see also \cite{groeger}), while expositions on \eqref{deGiorgiNash} can be found in the textbooks \cite{giltrud} or \cite{laduellipt}. \\

{\bf Statement of the main result:} We next formulate our main results, where we adopt the following way of speaking: We say that $\sigma$ is subject to (A) if it satisfies (A1), (A2) or (A2'), and (A3), depending on the context. The first result shows the existence of integrable second derivatives.
\begin{theo}\label{main}
	We assume that $\kappa \in C^{0,1}(\Omega; \, \mathbb{R}^{d\times d})$ satisfies the uniform ellipticity and boundedness assumptions \eqref{ellipt}, and that $\sigma$ is subject to (A). Moreover, for some $p \in ]\max\{d/(d-1), \, p_0^{\prime}\}, \, p_0[$ with $p_0$ according to condition \eqref{Meyers}, we assume that:
	\begin{enumerate}[(a)]
\item $f \in L^{p^{\rm \#}}(\Omega)$, and the function $E_1f: \, x \mapsto |x| \, f(x)$ belongs to $L^{p}(\Omega)$;
\item $Q  \in B^{p^{\rm \#}_S,p^{\rm \#}_S}_{1-1/p^{\rm \#}_S}(\Omega)$, and with $\beta_0 = \min\{1,(d-1)/p\}$ the function $E_{\beta_0}Q: \, x \mapsto |x|^{\beta_0} \, Q(x)$ belongs to $B^{p,p}_{1-1/p}(S)$.
\end{enumerate}
If $u$ is a local weak solution to \eqref{P1}, \eqref{P2}, then it belongs to $W^{1,p}_{\text{loc}}(\Omega)$ and to $W^{2,q}_{\text{loc}}(\Omega \setminus S)$ for all $1 \leq q <p^{\rm \#}$. Whenever $\Omega^{\prime\prime} \subset\!\!\subset \Omega^{\prime} \subset \!\! \subset \Omega$ (compact inclusions), we have
	\begin{align*}
	\|D^2 u\|_{L^q(\Omega^{\prime\prime} \setminus S)}  \leq C \, ( & \|f\|_{L^{p^{\rm \#}}(\Omega)} +\|E_1f\|_{L^p(\Omega)} + \|Q\|_{B^{p^{\rm \#}_S,p^{\rm \#}_S}_{p^{\rm \#}_S}(S)} + \|E_{\beta_0}Q\|_{B^{p,p}_{1-1/p}(S)} \\[-0.2cm]
	& + \|u\|_{L^p(\Omega^{\prime})})\, ,
		\end{align*}
	with $C$ depending on $q,p$ on $\beta_0$ and on ${\rm dist}(\Omega^{\prime},\partial \Omega)$, ${\rm dist}(\Omega^{\prime\prime},\Omega^{\prime})$.
\end{theo}
The next main result is the higher regularity of $\nabla u$ and $D^2u$ if the data $f$ and $Q$ are sufficiently regular. The argument relies on exploiting the H\"older continuity of the solution and an interpolation result of Miranda-Nirenberg type (see \cite{nirenberg66b}). In the cusp case, the domain $\Omega_2$ does not satisfy the cone condition, hence we rely on an inequality proved by own means in Prop.\ \ref{GNcusp}. Depending on the steepness of the cusp, this method might lead to a regularity upgrade. For simplicity, we restrict to the case $Q = 0$. 
\begin{theo}\label{main2}
We assume that $\sigma$ satisfies (A) and that $u$ is a local weak solution to \eqref{P1}, \eqref{P2}. Assume that $f \in L^{s_0}(\Omega)$ with $s_0 > d/2$, and that there are $s_1 > p_0^{\prime}$, $s_1 \geq s_0$ and $0 \leq \beta_1 \leq 1$ such that $\|E_{\beta_1} f\|_{L^{s_1}(\Omega)} < + \infty$.
We let $m_0 := \min\{p_0,s_1\}$ and $\lambda_0 := \min\{\alpha_0, 2-d/s_0\}$ with the constant $0 < \alpha_0\leq 1$ of condition \eqref{deGiorgiNash}. With $0 < \theta  \leq 1$ from (A3), we define $\theta_1 := 1$ and $\theta_2 := \theta$ and
\begin{align*}
	r_i := 
	\begin{cases}  d \, \min\Big\{\Big(\dfrac{m_0}{d+(\beta_1-1)\, m_0}\Big)^{\circ}, \, \Big(\dfrac{1+(\theta_i-1)/d}{1-\lambda_0}\Big)^{\circ}\Big\} & \text{ for } m_0 < + \infty \, ,\\
	+ \infty & \text{ for } m_0 = + \infty \, .
	\end{cases} \quad ( i = 1,2)\, .
\end{align*}
where for $x > 0$ and $y \in \mathbb{R}$, we set $(x/y)^{\circ} := x/y$ for $y >0$, $(x/y)^\circ$ arb.\ large positive for $y =0$, $(x/y)^{\circ} := +\infty$ for $y < 0$. Then $\nabla u \in L^r_{\text{loc}}(\Omega_i; \, \mathbb{R}^d)$ for all $1 \leq r < r_i$ and $D^2 u \in L^{q}_{\text{loc}}(\Omega_i; \, \mathbb{R}^{d\times d})$ for all $1 \leq q < r_id/(r_i+d)$. For all $\Omega^{\prime\prime}\subset\!\!\subset \Omega^{\prime}\subset\!\!\subset\Omega$, there is a constant $c = c(r_i,q,\Omega^{\prime\prime},\Omega^{\prime},k_0,k_1,s_0,s_1,\beta_1)$ such that
\begin{align*}
&	\|D^2 u\|_{L^{q}(\Omega^{\prime\prime} \cap \Omega_i)} + 	\|\nabla u \|_{L^{r}(\Omega^{\prime\prime}\cap \Omega_i)}  \leq c \, (\|E_{\beta_1} f\|_{L^{s_1}(\Omega)} +\|f\|_{L^{s_0}(\Omega)}+ \|u\|_{L^{z}(\Omega^{\prime})}) \, ,
\end{align*}
where we can take $z = (d \, s_0/(d-s_0))^{\circ}$.
\end{theo}
\begin{rem}
	For $f = 0$ and $Q = 0$, we obtain the following regularity, involving the constant $\alpha_0$ of \eqref{deGiorgiNash} and $\theta$ of (A3):
	\begin{enumerate}
		\item $\nabla u \in L^r(\Omega_1)$ for all $r < \min\{p_0^*, d/(1-\alpha_0)\}$ and $D^2u \in L^q(\Omega_1)$ for all $q <\min\{p_0,\, d/(2-\alpha_0)\}$. 
		\item $\nabla u \in L^r(\Omega_2)$ for all $r < \min\{p_0^*, (d+\theta-1)/(1-\alpha_0)\}$ and $D^2u \in L^q(\Omega_2)$ for all $q < \min\{p_0,\, (d+\theta-1)/[2-\alpha_0-(\theta-1)/d]\}$. 
	\end{enumerate}
Hence, we obtain $\nabla u$ integrable to a power larger than $d$ in $\Omega_1$, while the same is correct in $\Omega_2$ only if the cusp is not too steep, viz.\ for $\alpha_0> (1-\theta)/d$. If $\theta$ is too small, the method even might give no improvement of regularity for the gradient in $\Omega_2$. For the physical dimensions $d = 2,3$ we reach $(d+\theta-1)/(1-\alpha_0) \leq p_0$ if $0< \theta \leq 1 + p_0 - d  - p_0 \, \alpha_0$. Thus, there are cusps where the method leads to no gain in regularity if the interval is not empty, which in 2D is the case if $\alpha_0 < 1 - 1/p_0$ and in 3D if $\alpha_0$ and $p_0-2$ are sufficiently small.
\end{rem}

\begin{rem}
	The results are compatible with the known counter-examples of \cite{meyers}, \cite{reh} in dimension $d= 2$, showing for polygonal interfaces that the constant $p_0$ might turn arbitrarily close to $2$, depending on the eigenvalues of $\kappa$ on both sides of the interface. However, our result also show for $d > 2$ that the solutions to interface problems might behave better under point singularities.   
\end{rem}

\section{Directional regularity}\label{direction}

We recall that $x = (\bar{x},  \, x_d)$ with $\bar{x} = (x_1,\ldots,x_{d-1})$ and we define $r = r(\bar{x}) := |\bar{x}|$. In order to make use of the properties (A) for $\sigma$, we restrict the choice of $x$ to the cylinder 
\begin{align}\label{Zylr0}
	\Omega_{R_0} := \{x \in \mathbb{R}^d \, : \, |\bar{x}| < R_0, \, |x_d| < \sigma(R_0)\} \, .
\end{align}
In the spirit of local regularity, we replace from now the domain $\Omega$ by its intersection with $\Omega_{R_0}$, which simplifies proofs and in particular allows to make use of the properties (A) for the entire $S \cap \Omega$.

Under these conditions, a unit normal field, continuous on $S\setminus\{0\}$, is given by
\begin{align}\label{nus}
	\nu^S(x) := \frac{1}{(1+(\sigma^{\prime}(r))^2)^{\frac{1}{2}}} \, \begin{pmatrix}
		- \sigma^{\prime}(r) \, \bar{x}/r\\
		1
	\end{pmatrix}\quad \text{ for } \quad x \in S \, .
\end{align}
We define a unit tangent field via
\begin{align}\label{taus}
	\tau^S(x) := \frac{1}{(1+(\sigma^{\prime}(r))^2)^{\frac{1}{2}}} \, \begin{pmatrix}
		\bar{x}/r\\
		\sigma^{\prime}(r) \, 
	\end{pmatrix} \quad \text{ for } \quad x \in S \, .
\end{align}
For a cusp, see (A2), the continuity of $\tau^S(x)$ over $S$ up to $x = 0$ is easily verified (but it fails for a cone!). Our next purpose is to construct extensions into $\Omega$ for $\tau^S$. To this aim we first introduce a few auxiliary functions.
\begin{lemma}\label{Mandmu}
	Let $\sigma$ be subject to (A). For $z \in ]-\sigma(R_0), \, \sigma(R_0)[ =: I_0$ we define $M(z) :=  \int_0^{|z|}\sigma^{\prime}(\sigma^{-1}(s)) \, \sigma^{-1}(s) \, ds$. Then $M$ belongs to $C^{1,1}(I_0) \cap C^2(I_0\setminus \{0\})$, is a strictly convex function with minimum at $0$, which satisfies $M(z) = M(-z)$ and $\theta \, z^2 \leq 2 \, M(z) \leq z^2$ for all $z \in I_0$.
For all $\rho \in ]0,R_0[$, we denote by $M^{-1}_+(\rho)$ the positive root of the equation $M(z) = \rho$. Then, the function $ M^{-1}_+$ is continuously differentiable and increasing in $]0,R_0[$, and it satisfies
\begin{align*}
	\sqrt{2\rho} \leq M^{-1}_+(\rho) \leq \sqrt{\frac{2\rho}{\theta}} \quad \text{ for all } \quad \rho \in ]0,R_0[ \, .
\end{align*}
\end{lemma}
\begin{proof}
	Since $\sigma^{\prime} >0$ in $]0,\, R_0]$, the function $\sigma$ is invertible over its domain.
	For $z \in  I_0$ we define 	\begin{align}\label{defofmu} 
		\mu(z) := {\rm sign}(z)\,	\sigma^{\prime}\Big(\sigma^{-1}(|z|)\Big) \, \sigma^{-1}(|z|) \, .
	\end{align}
	Owing to the assumption (A3), we can verify\footnote{(A3) yields the differential inequality $(\ln \sigma^{\prime}(\rho))^{\prime} \geq (\theta-1) /\rho$, which we can integrate on each interval $[\rho, \, R_0]$.} that $\sigma^{\prime}(\rho) \leq \sigma^{\prime}(R_0) \, (R_0/\rho)^{1-\theta}$ for $0< \rho\leq R_0$. Hence $\rho \, \sigma^{\prime}(\rho)$ tends to zero for $\rho \to 0$. Since $\sigma^{-1}(|z|) \to 0$ for $z \to 0$, we see that $\mu$ is bounded and continuous in zero with $\mu(0) = 0$.  Moreover, we compute that
	\begin{align*}
		\mu^{\prime}(z) = 
		\Big(	1 + \frac{\sigma^{\prime\prime}(\rho)\, \rho}{\sigma^{\prime}(\rho)}\Big)_{\rho = \sigma^{-1}(|z|)} \text{ for } |z| > 0 \, .
	\end{align*}	
	In view of (A3) again, we see that $\mu^{\prime}$ is bounded on $\Omega$ with $\theta \leq \mu^{\prime} < 1$, hence also
	\begin{align}\label{mulineargrowth}
		\frac{\mu(z)}{z} \geq \theta > 0 \quad \text{ for  all } \quad z \in [-\sigma(R_0), \, \sigma(R_0)] \, .
	\end{align}
Since $M$ is the primitive of $\mu$ that vanishes in zero, we easily verify that $M(z) = M(-z)$, and as $M^{\prime\prime} = \mu^{\prime} \geq \theta$ it follows that $M$ is strictly convex. By the Taylor formula, $M(z) = M^{\prime\prime}(\lambda \, z) \, z^2/2$ with a $\lambda \in ]0,1]$, hence $\theta \, z^2\leq 2M(z) \leq z^2$ as claimed.

Since $M$ is strictly convex and even, the equation $M(z) = \rho$ possesses for $\rho \in ]0,R_0[$ exactly two solutions $\pm z$. We denote $M^{-1}_+(\rho)$ the positive root, which is unique in a neigbourhood of $\rho$. Since $M(M^{-1}_+(\rho)) = \rho$, the implict function theorem guarantees that $ M^{-1}_+$ is differentiable with $(M^{-1}_+)^{\prime}(\rho) = 1/M^{\prime}(M^{-1}_+(\rho)) = 1/\mu(M^{-1}_+(\rho))$. Since $\mu$ is positive on $]0,\sigma(R_0)]$ this shows that $M^{-1}_+$ is increasing.

Finally, since we have shown that $\theta \, z^2 \leq 2M(z) \leq z^2$, we also have for $\rho \in [0, \, R_0]$ arbitrary that $\sqrt{2\rho} \leq M^{-1}_+(\rho) \leq \sqrt{2\rho/\theta}$ as claimed.
\end{proof}
As an illustration, note that for the case $\sigma(\rho) = \rho^{\theta}$, we easily verify that $M(z) = \theta\, |z|^2/2$.
\begin{lemma}\label{tauregu}
	Let $\sigma$ be subject to the same conditions as in Lemma \ref{Mandmu}.
Then there is a vector field $\tau \in \bigcap_{1\leq p < d} W^{1,p}(\Omega; \, \mathbb{R}^d)$ such that $\tau(x) = \tau^S(x)$ for all $x \in S$. For all $1 \leq p < d$, there is $c = c(p,R_0) > 0$ such that $\|\nabla \tau\|_{L^{p}(\Omega)} \leq c \, \theta^{-1}$. Moreover we have $\tau \in C^1(\Omega \setminus\{0\}; \, \mathbb{R}^d)$.
\end{lemma}
\begin{proof}
We define
\begin{align}\label{gedefin}
	g(x) := \sqrt{\frac{1}{2} \, |\bar{x}|^2 + M(x_d)}\quad \text{ for } \quad x \in \Omega \, ,
\end{align}
and Lemma \ref{Mandmu} implies that $g(x) = g(-x)$ and 
\begin{align}\label{gradientofg}
	\frac{\nabla g(x)}{|\nabla g(x)|} = \frac{1}{(|\bar{x}|^2 + \mu^2(x_d))^{\frac{1}{2}}} \, \begin{pmatrix}
		\bar{x}\\
		\mu(x_d)
	\end{pmatrix} \quad \text{ for } \quad x\neq 0 \, .
\end{align}
For $x \in S$, we have $x_d = \sigma(|\bar{x}|)$, and we get 
\begin{align*}
	\frac{\nabla g(x)}{|\nabla g(x)|} = \frac{1}{(r^2 + (\sigma^{\prime}(r) \, r)^2)^{\frac{1}{2}}} \, \begin{pmatrix}
		\bar{x}\\
		\sigma^{\prime}(r) \, r
	\end{pmatrix}  = \tau^S(x) \quad \text{for} \quad x \in S \, .
\end{align*}
Thus, defining $\tau(x) = \nabla g(x)/|\nabla g(x)|$ for $x \in \Omega \setminus S$, we find a normalised extension of the vector field $\tau^S$. The generalised derivatives $\{\partial_{x_i}\tau_k\}$ form a symmetric matrix given as
\begin{gather*}
\partial_{x_i}\tau_k(x) = \frac{1}{\sqrt{r^2 + \mu^2(x_d)}} \, \Big(\delta_{ik} - \frac{x_ix_k}{r^2 + \mu^2(x_d)}\Big) \quad \text{ for } \quad i,k=1,\ldots,d-1\\
\partial_{x_i}\tau_d = \partial_{x_d}\tau_i = - \frac{ x_i \, \mu(x_d)\, \mu^{\prime}(x_d)}{(r^2 + \mu^2(x_d))^{\frac{3}{2}}}  \quad \text{ for } \quad i =1,2\,, \qquad
\partial_{x_d}\tau_d = \frac{ \mu^{\prime}(x_d) \, r^2}{(r^2 + \mu^2(x_d))^{\frac{3}{2}}}\, .
\end{gather*}
Thus, it is readily proven that
\begin{align*}
|\partial_{x_i}\tau_k(x)| \leq \frac{\max\{1,\, \|\mu^{\prime}\|_{L^{\infty}}\}}{\sqrt{r^2 + \mu^2(x_d)}} \leq \frac{1}{\sqrt{r^2 + \mu^2(x_d)}} \, 
\end{align*}
Since \eqref{mulineargrowth} implies that $\mu^2(x_d) >\theta^2 \, x_d^2 $, it holds that
\begin{align*}
	\sqrt{r^2 + \mu^2(x_d)} \geq \theta \, |x| \quad \text{ in } \quad \Omega \, .
\end{align*}
Therefore, we see that $\|\nabla \tau\|_{L^p} \leq \theta^{-1} \, \|\frac{1}{|\cdot|}\|_{L^p}$. Hence $\tau \in W^{1,p}(\Omega; \, \mathbb{R}^d)$ for all $1\leq p < d$, and the estimate easily follows. Obviously, $\tau$ is continuously differentiable outside of zero.
\end{proof}
We will prove regularity for data in certain weighted classes. We recall that for a function $f$ defined in $\Omega$, and for $\beta >0$, we let $E_{\beta}f$ be the function $x \mapsto f(x) \, |x|^{\beta}$.
We first consider $Q = 0$ and employ these considerations to prove the following statement.
\begin{prop}\label{v1prop}
	With a constant $p_0$ occurring in \eqref{Meyers}, we consider an arbitrary number $r \in ]\max\{p_0^{\prime},d/(d-1)\}, \, p_0^*[$. Assume that $f$ satisfies $\|E_{\beta_1}f\|_{L^{s_1}(\Omega)} < + \infty$ with parameters $0 \leq \beta_1 \leq 1$ and $\max\{p_0^{\prime}, \, r^{\rm \#}\} < s_1 \leq + \infty$ subject to the restriction
	\begin{align*}
		\begin{cases}
	\beta_1 \leq 1 & \text{ for } r \leq \min\{p_0,s_1\}\, ,\\
		\beta_1 < 1 + d/r - d/\min\{p_0,s_1\} & \text{ otherwise} \, .
		\end{cases}
	\end{align*}	
	Suppose that $u \in W^{1,r}_{\text{loc}}(\Omega)$ is a local weak solution to \eqref{P1}, \eqref{P2}. Then the function $v_1 := \tau \cdot \nabla u$ belongs to $\bigcap_{1\leq q < \frac{dr}{d+r}} W^{1,q}_{\text{loc}}(\Omega)$. For each $ 1\leq q < dr/(d+r)$ and $\Omega^{\prime\prime} \subset\!\!\subset \Omega^{\prime} \subset\!\!\subset \Omega$ there is $c = c(r,q,s_1,\beta_1,\Omega^{\prime},k_0,k_1)$ such that
	\begin{align*}
		\|\nabla v_1\|_{L^{q}(\Omega^{\prime\prime})} \leq c \, (\|\nabla u\|_{L^{r}(\Omega^{\prime})} + \|E_{\beta_1} f\|_{L^{s_1}(\Omega^{\prime})}) \, .
	\end{align*} 
\end{prop}
\begin{proof}
	As a preliminary, we construct parallel modifications of the vector field $\tau$ which temper the singularity. For $\epsilon \geq 0$, and $t \geq 0$, we define 
	\begin{align*}
		f_{\epsilon}(t) = \begin{cases}
			0 & \text{ for } 0 \leq t \leq \epsilon\\
			t-\epsilon & \text{ for } \epsilon < t < +\infty
		\end{cases}
		\, ,
	\end{align*}
	and for fixed $0 < \beta \leq 1$ we let $\tau^{\epsilon}(x) := f_{\epsilon}(|x|^{\beta}) \, \tau(x)$. Using the representation of $\nabla \tau$ obtained in Lemma \ref{tauregu}, in connection with the identity $\nabla \tau^{\epsilon}(x) = f_{\epsilon}(|x|^{\beta}) \, \nabla \tau + \beta \, |x|^{\beta-1} \, f_{\epsilon}^{\prime}(|x|^{\beta}) \, \frac{x}{|x|} \otimes \tau(x)$, it is readily shown that 
	\begin{align}\label{estimtaueps}
		|\nabla \tau^{\epsilon}(x)| \leq \Big(\frac{1}{\theta} + \beta\Big)\, \chi_{\{y \,  : |y|^{\beta} > \epsilon\}}(x)  \, |x|^{\beta-1} \quad \text{ for all } \quad x \in \Omega, \, \epsilon > 0 \,.
	\end{align}
	Here $\chi_A$ denotes the characteristic function of a set $A$.
	We remark that $\tau^{\epsilon}$ is of class $W^{1,\infty}(\Omega;\, \mathbb{R}^d)$ for all $\epsilon >0$, since $|\nabla \tau^{\epsilon}(x)| \leq c \, \epsilon^{\beta-1}$.\\

	For $\epsilon > 0$, the interface $S\setminus B_{\epsilon}(0)$ is of class $\mathscr{C}^{\infty}$, and it is therefore classical to show that every weak solution to \eqref{P1}, \eqref{P2} belongs to $W^{2,2}_{\text{loc}}(\Omega \setminus B_{\epsilon}(0))$ (Ref.\ \cite{stamptransmi}). Using that $\tau^{\epsilon}$ belongs to $W^{1,\infty}(\Omega; \, \mathbb{R}^d)$ and is supported in $\Omega \setminus B_{\epsilon^{1/\beta}}(0)$, we then obtain that the function $v_{\epsilon} := \tau^{\epsilon} \cdot \nabla u$ belongs to $W^{1,2}_{\text{loc}}(\Omega)$.
	
	Moreover, for $\phi \in C^{\infty}_c(\Omega)$ arbitrary, we can test the equation \eqref{P1} with functions of the form $\tau^{\epsilon}\cdot \nabla \phi$. Note that $\tau^{\epsilon}$ is parallel to $\tau$, and therefore it is tangential on $S \setminus \{0\}$. This allows to show by means of an integration by parts and standard calculus that $v_{\epsilon}$ satisfies (see \cite{druetboheme}, \cite{druet2018regularity} for similar computations)
	\begin{align}\label{v1distribution}
		\int_{\Omega} \kappa \nabla v_{\epsilon} \cdot \nabla \phi  \, dx = \int_{\Omega}\{ D_{\kappa}(\tau^{\epsilon}) \nabla u - (\tau^{\epsilon}\cdot \nabla )\kappa\,  \nabla u - f \, \tau^{\epsilon}\} \cdot \nabla \phi \, dx \quad \text{ for all } \quad \phi \in C^{\infty}_c(\Omega)\, ,
	\end{align}
	where
	$D_{\kappa}(\tau^{\epsilon}) = \kappa \, \nabla \tau^{\epsilon} +  (\nabla \tau^{\epsilon})^{\sf T} \, \kappa - \divv (\tau^{\epsilon}) \, \kappa$. Invoking the property \eqref{Meyers}, for $\Omega^{\prime\prime} \subset\!\!\subset \Omega^{\prime} \subset\!\!\subset \Omega$ and for $p_0^{\prime} < s < p_0$ arbitrary, we obtain the estimate
	\begin{align*}
		\|\nabla v_{\epsilon}\|_{L^{s}(\Omega^{\prime\prime})} \leq c \, (\|D_{\kappa}(\tau^{\epsilon}) \nabla u\|_{L^{s}} + \|f \, \tau^{\epsilon}\|_{L^{s}} + \||\tau^{\epsilon}| \, \nabla u\|_{L^{s}}) \, ,
	\end{align*}
	where we use that $\kappa$ is of class $C^{0,1}(\Omega\setminus S;\mathbb{R}^{d\times d})$. Here $L^s$ stands for $L^{s}(\Omega^{\prime})$. 
	Invoking the preliminary consideration, we next have $|D_{\kappa}(\tau^{\epsilon})(x)| \leq c \, k_1 \, |x|^{\beta-1}$. 
	We restrict the choice of $s$ and $\beta$ by the conditions
	\begin{align}\label{CON1}
	\beta > 1 - d \, \Big(\frac{1}{p_0^{\prime}}-\frac{1}{r}\Big) \quad \text{ and } \quad p_0^{\prime} < s < \frac{dr}{d + (1-\beta) \, r} \, 
	\end{align}
so that we can rely on the properties $s \leq r$ and $(\beta-1)rs/(r-s) > - d$. By means of H\"older's inequality, we infer that
	\begin{align*}
		\|D_{\kappa}(\tau^{\epsilon}) \nabla u\|_{L^s} & \leq c \, k_1 \,\||\cdot|^{\beta-1} \, \nabla u\|_{L^s} \\
		& \leq c \, k_1\,  \|\nabla u\|_{L^r} \, \Big(\int_{\Omega} |x|^{(\beta-1)rs/(r-s)} \, dx\Big)^{\frac{1}{s}-\frac{1}{r}} \leq C \, k_1 \, \|\nabla u\|_{L^r} \, ,
	\end{align*}
	which implies that $\|\nabla v_{\epsilon}\|_{L^{s}(\Omega^{\prime\prime})} \leq c \, (\|\nabla u\|_{L^{r}} + \|E_{\beta}f\|_{L^{s}})$. If we further require that
	\begin{align}\label{CON2}
	\beta \geq \beta_1 \quad \text{ and }  \quad s \leq s_1 \, , 
	\end{align}
then the norm of $f$ can be controlled via
	\begin{align*}
		\|E_{\beta}f\|_{L^{s}} \leq \sup_{x \in \Omega} |x|^{\beta-\beta_1} \, |\Omega|^{1-\frac{s}{s_1}} \, \|E_{\beta_1} f\|_{L^{s_1}} \, ,
	\end{align*}
from which we obtain that $
	\|\nabla v_{\epsilon}\|_{L^{s}(\Omega^{\prime\prime})} \leq c \, (\|\nabla u\|_{L^{r}} + \|E_{\beta_1}f\|_{L^{s_1}})$. Letting $\epsilon \to 0$, we obtain for $w := E_{\beta} v_1$ the bound 
\begin{align}\label{ESTIGROUND}
	\|\nabla w\|_{L^{s}(\Omega^{\prime\prime})} \leq c \, (\|\nabla u\|_{L^{r}} + \|E_{\beta_1}f\|_{L^{s_1}}) \, .
\end{align}
We let $m_0 := \min\{p_0,s_1\}$, and since we assume that $s_1 > p_0^{\prime}$ and $r > p_0^{\prime}$, we observe that $m_0> p_0^{\prime}$. Since we also assume that $r< p_0^*$ and $s_1 > r^{\rm \#}$, we have $r < m_0^{\rm \#}$, hence $1 + d/r -d/m_0 > 0$. Next we choose arbitrary $\beta$ and $s$ as follows:
\begin{align}\label{CON3}
&	\begin{cases}
		\beta_1 \leq \beta \leq 1 & \text{ if }  m_0 \geq r \text{ and }  \beta_1 > 1 - d \, (\frac{1}{p_0^{\prime}}-\frac{1}{r}) \, ,\\
		1 - d \, (\frac{1}{p_0^{\prime}}-\frac{1}{r}) < \beta \leq 1 & \text{ if } m_0 \geq r  \text{ and }  \beta_1 \leq  1 - d \, (\frac{1}{p_0^{\prime}}-\frac{1}{r})\, ,\\
\max\{\beta_1,\,  1 - d \,(\frac{1}{p_0^{\prime}}-\frac{1}{r})\} < \beta <  1 + \frac{d}{r} - \frac{d}{m_0} & \text{ if } m_0 < r \, ,\\
\end{cases}\\
\label{CON3bis} & p_0^{\prime} < s < \frac{dr}{d + (1-\beta) \, r} \, . 	
	\end{align}
We then can verify that \eqref{CON1} and \eqref{CON2} are valid. In particular, if $q < rd/(r+d)$, we can always find $s>q$ such that $\beta \, qs/(s-q)< d$. Indeed, the latter condition is equivalent to $q < s/(1+\beta s/d)$, and the threshold of this function for $s \nearrow dr/[d+(1-\beta) \, r]$ is the number $dr/(d+r)$ which is strictly larger than $q$. For given $q$, we next specialise to such $s$ satisfying \eqref{CON3bis} and $\beta \, qs/(s-q)< d$.

We then observe that $\nabla v_1 = |x|^{-\beta} \, \nabla w - \beta \, |x|^{-1} \, (\tau \cdot \nabla u) \, \nabla |x|$. Hence, H\"older's inequality implies that
	\begin{align*}
		\|\nabla v_1\|_{L^q(\Omega^{\prime\prime})} \leq & \beta \, \||x|^{-1} \, \nabla u\|_{L^q(\Omega^{\prime\prime})} + \| |x|^{-\beta} \, \nabla w\|_{L^q(\Omega^{\prime\prime})} \\
		\leq & \beta \,\Big(\int_{\Omega} |x|^{-rq/(r-q)} \, dx\Big)^{\frac{1}{q}-\frac{1}{r}} \, \|\nabla u\|_{L^r} + \Big(\int_{\Omega} |x|^{-\beta qs/(s-q)} \, dx\Big)^{\frac{1}{s}-\frac{1}{q}} \, \|\nabla w\|_{L^s(\Omega^{\prime\prime})}\\
		\leq & C \, (\|\nabla u\|_{L^r} +\|E_{\beta_1 }f\|_{L^{s_1}} )\, . \qedhere
	\end{align*}
\end{proof}
We note that the condition $f \in L^{p^{\rm\#}}(\Omega)$ ensures that the functional $\phi \mapsto \int_{\Omega} f \, \phi \, dx$ belongs to $[W^{1,p^{\prime}}(\Omega)]^*$. With respect to the problem \eqref{P1}, \eqref{P2}, where we assume first $Q = 0$, the result \eqref{Meyers} implies that if $p^{\prime}_0 < p < p_0$ then every local weak solution $u$ belongs to $W^{1,p}_{\text{loc}}(\Omega)$ and satisfies the estimate
\begin{align}\label{ugroeg}
	\|\nabla u\|_{L^{p}(\Omega^{\prime\prime})} \leq c(p,\Omega^{\prime\prime}) \, (\|f\|_{L^{p^{\rm\#}}(\Omega^{\prime})}+\|u\|_{L^p(\Omega^{\prime})}) \quad \text{ for all } \quad \Omega^{\prime\prime} \subset\!\!\subset \Omega^{\prime} \subset\!\!\subset  \Omega \,.
\end{align}
Still considering $Q = 0$ we have the following consequence.
\begin{coro}\label{v1coro}
Let $\max\{p_0^{\prime},d/(d-1)\} < p < p_0 $ with the constant $p_0$ occurring in \eqref{Meyers}. We assume that $f \in L^{p^{\rm\#}}(\Omega)$ with $E_1f \in L^{p}(\Omega)$ and that $u$ is a local weak solution to \eqref{P1}, \eqref{P2}, satisfying $u \in L^{p}_{\text{loc}}(\Omega)$. Then the function $v_1 := \tau \cdot \nabla u$ belongs to $\bigcap_{1\leq q < p^{\rm \#}} W^{1,q}_{\text{loc}}(\Omega)$. For each $ 1\leq q < p^{\rm \#}$ and $\Omega^{\prime\prime} \subset\!\!\subset \Omega^{\prime} \subset\!\!\subset  \Omega$ with compact inclusions, there is $c = c(p,q,\Omega^{\prime\prime},k_0,k_1)$ such that
\begin{align*}
		\|\nabla v_1\|_{L^{q}(\Omega^{\prime\prime})} \leq c \, (\|f\|_{L^{p^{\rm\#}}(\Omega^{\prime})} + \|E_1f\|_{L^{p} (\Omega^{\prime})}+\|u\|_{L^p(\Omega^{\prime})}) \, .
\end{align*} 
\end{coro}
\begin{proof}
We apply the Proposition \ref{v1prop} with $r = p $, $s_1 = p$ and $\beta_1 = 1$. For any $1 \leq q < p^{\rm \#}$ we obtain the inequality $\|\nabla v_1\|_{L^q(\Omega^{\prime\prime})} 	\leq C \, (\|\nabla u\|_{L^p} +\|E_{1}f\|_{L^p} )$. Combining it with \eqref{ugroeg}, the claim follows.
\end{proof}

\section{Dimension reduction}\label{dmoinsun}

The next idea is to obtain Sobolev--space-estimates for the derivatives of $u$ in directions orthogonal to the vector $\tau$. To this aim, for $R > 0$, we denote by $\Gamma_R$ the $R-$level-set of the function $g$ introduced in \eqref{gedefin}, that is,
\begin{align}\label{levelset}
	\Gamma_R := \Big\{x \in \Omega \, : \, \frac{1}{2} \, |\bar{x}|^2+ M(x_d) = R^2\Big\}\, .
\end{align}
Recall that we have identified $\Omega$ with $\Omega_{R_0}$ of \eqref{Zylr0}. A unit normal for $\Gamma_R$ is obviously given by $ \nabla g/|\nabla g| = \tau$. 
Next we state a few useful properties that further characterise $\Gamma_R$. To this aim we also define a strictly positive, nonincreasing function
	\begin{align}\label{phinull}
		\phi_0(R) := \inf_{\rho \in ]0,R[} \frac{\sigma(\rho)}{\rho} \quad \text{ for }\quad 0 < R < R_0 \, .
		\end{align}
Observe that $\lim_{R\rightarrow 0+}\phi_0(R) = + \infty$ for (A2), while $\lim_{R\rightarrow 0+}\phi_0(R) = \sigma_0^{\prime} >0$ for (A2') due to standard rules of calculus.  

\begin{minipage}{\textwidth}
	\centering
		\includegraphics[width = 0.6\textwidth]{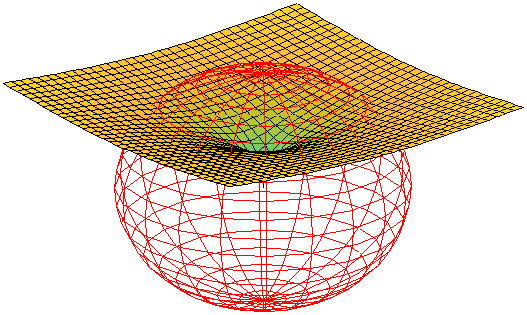}
	\captionof{figure}{The ellipsoid $\Gamma_R$ represented in red is intersecting the graph of $\sigma(r) = r^{\frac{1}{3}}$ transversally, here for the value $R = 5/(2\sqrt{2})$.}
\end{minipage}

\begin{lemma}\label{curve}
	We assume that $\sigma$ satisfies (A). For all $R \in  ]0,R_0]$, the surface $\Gamma_R$ is a regular closed manifold of class $\mathscr{C}^2$. Moreover, $\Gamma_R$ is the boundary of the convex set
	\begin{align}\label{Cr}
C_R := \Big\{x \in \Omega \, : \, \frac{1}{2}\, |\bar{x}|^2 + M(x_d) < R^2\Big\} \ni \{0\} \, .
\end{align}
The intersection $\Gamma_R \cap S$ is a $(d-2)-$dimensional sphere characterised as $\{x \in \Omega \, : \, |\bar{x}| = \hat{\rho}(R), \, x_d = \hat{z}(R)\}$ with two positive functions $\hat{\rho}, \, \hat{z}$. In particular, there is $a_0 >0$ such that
	\begin{align*}
	a_0 \, R^{\frac{1}{\theta}} \leq \hat{\rho}(R) \leq \frac{\sqrt{2}R}{(1+\theta \, (\phi_0(\sqrt{2}R))^2)^{\frac{1}{2}}} \quad \text{ for all } \quad 0 < R < \frac{R_0}{\sqrt{2}} \, .
\end{align*}
\end{lemma}
 \begin{proof}
Starting with the definition \eqref{gedefin} we compute that
\begin{align*}
	|\nabla g(x)| = \frac{1}{2} \, \sqrt{\frac{|\bar{x}|^2+ \mu(x_d)^2}{\frac{1}{2} \, |\bar{x}|^2+ M(x_d)}} \, .
\end{align*}
Exploiting Lemma \ref{Mandmu}, in particular \eqref{mulineargrowth}, we have $\mu(x_d)^2 \geq \theta^2 \, x_d^2 \geq 2 \, \theta^2 \, M(x_d)$ and $ \mu(x_d)^2 \leq \max |\mu^{\prime}|  \, x_d^2 \leq 2 \, M(x_d)/\theta^2$. This implies that
\begin{align}\label{nablagbounds}
	\frac{\theta}{\sqrt{2}} \leq |\nabla g(x)| \leq \frac{1}{\sqrt{2}\, \theta} \quad \text{ for all } \quad x \in \Omega \, .
\end{align}

Since $\nabla g(x) \neq 0$ for $x \in \Omega$, the level sets of $g$ are regular manifolds of class $\mathscr{C}^2$. Recall that $M$ is convex (cf.\ Lemma \ref{Mandmu}). Thus the interior of $\Gamma_{R}$ is convex, and due to the fact that $g(x) = g(-x)$, the sets $C_R$ are convex sets containing zero.  

A point $x$ belongs to $S \cap \Gamma_R$ if and only if it obeys the equations $x_d = \sigma(r)$ and $r^2/2 + M(x_d) = R^2$. Thus $(\sigma^{-1}(x_d))^2/2 + M(x_d) = R^2$ and $r^2/2 + M(\sigma(r)) = R^2$.

Since $M$, $\sigma$ and $\sigma^{-1}$ are increasing and continuously differentiable on positive arguments, these equations determine uniquely $x_d = \hat{z}(R)$ and $|\bar{x}| = \hat{\rho}(R)$ with differentiable functions.
Employing Lemma \ref{Mandmu}, we have
\begin{align*}
	R^2 = \frac{1}{2}\hat{\rho}^2 + M (\sigma(\hat{\rho})) \leq \frac{1}{2} \, (\hat{\rho}^2 + (\sigma(\hat{\rho}))^2) \, .
	\end{align*}
Using the assumption (A3), we estimate 
\begin{align*}
	\sigma(\hat{\rho}) \leq \frac{\sigma(R_0)}{\theta} \, \left(\frac{\hat{\rho}}{R_0}\right)^{\theta} \, ,
\end{align*}
and combining both inequalities we get
\begin{align*}
	R^2 \leq   \frac{1}{2R_0^{2\theta}} \, \Big(R_0^{2} +\frac{\sigma^2(R_0)}{\theta^{2}}\Big) \, \hat{\rho}^{2\theta} \, ,
\end{align*}
which shows that $\hat{\rho} \geq a_0 \, R^{1/\theta}$ with $a_0 > 0$ constant. Note that obviously $\hat{\rho}(R) \leq \sqrt{2}R$. We invoke Lemma \ref{Mandmu} again, and restricting to $R < R_0/\sqrt{2}$ this yields
\begin{align*}
	R^2 = & \frac{1}{2}\hat{\rho}^2 + M (\sigma(\hat{\rho})) \geq \frac{1}{2} \,  \, (\hat{\rho}^2+ \theta\,  (\sigma(\hat{\rho}))^2) =  \frac{1}{2} \, (1+\theta \, (\sigma(\hat{\rho}) /\hat{\rho})^2) \,  \hat{\rho}^2\\
		\geq &  \frac{1}{2} \,  \, (1 + \theta \, \phi_0(\sqrt{2}R))^2) \,  \hat{\rho}^2\, . \qedhere
\end{align*}
 \end{proof}
\begin{coro}\label{w1infty}
	Under the assumptions of Lemma \ref{curve}, 
	there is $\lambda_1 >0$ such that the vector fields $\tau$ and $\nu$  satisfy
	\begin{align*}
	 |\nu(x)-\nu(y)| \leq \lambda_1 \,\frac{|x-y|}{\hat{\rho}(R)} \quad \text{ and } \quad |\tau(x)-\tau(y)| \leq \lambda_1 \, \frac{|x-y|}{(1+|\sigma^{\prime}(\hat{\rho}(R))|^2)^{\frac{1}{2}}\, \hat{\rho}(R)}\, ,
	 \end{align*}
for all $x,\, y \in \Gamma_R\cap S$, $0< R \leq R_0$.
\end{coro}
\begin{proof}
We let $\partial_s := t \cdot \nabla$, where $t$ is an arbitrary unit tangent vector on the sphere $\Gamma_R \cap S$.
Using the representations \eqref{nus} and \eqref{taus}, we can compute that
\begin{align*}
	\partial_s \nu = \frac{-\sigma^{\prime}(r)}{(1+(\sigma^{\prime}(r))^2)^{\frac{1}{2}}} \, \frac{t}{r} ,\quad \partial_s \tau = \frac{1}{(1+(\sigma^{\prime}(r))^2)^{\frac{1}{2}}} \, \frac{t}{r}\, , 
\end{align*}
	and since $t$ possesses unit length and $r = \hat{\rho}(R)$ is constant on $\Gamma_R\cap S$, it follows that
	\begin{align*}
		|\partial_s \nu|\leq \frac{1}{\hat{\rho}(R)}  , \quad |\partial_s \tau|\leq \frac{1}{(1+(\sigma^{\prime}(\hat{\rho}(R)))^2)^{\frac{1}{2}}\, \hat{\rho}(R)} \, .
	\end{align*}
For $x,y \in \Gamma_R\cap S$ we denote by by $L(x,y)$ the geodesic distance between $x$ and $y$ in the sense of the chosen orientation $t$. Since $L(x,y) \leq \lambda_1 \, |x-y|$, the claim follows.	
\end{proof}
Before starting the next discussions, we need an auxiliary Lemma, which shall help us transforming estimates on the curved manifolds $\Gamma_R$ to bulk estimates. For $0 \leq R \leq R_0$ we introduce the sets
\begin{align*}
	S_R := S \cap C_R \, .
\end{align*}
\begin{lemma}\label{Twoformeln}
For all $f \in W^{1,1}(\Omega)$ the function $R \mapsto \int_{C_R} f   \, dx$ is differentiable in every point $R \in (0,R_0)$, and    
\begin{align*}
	\int_{\Gamma_R} f \, |\nabla g|^{-1} \, da = \frac{d}{dR} \, \int_{C_R} f   \, dx  \, .
\end{align*}
Moreover, if $f \in W^{1,r}(\Omega)$ with $r > 2$, then the function $R \mapsto \int_{S_R} f \, da$ belongs to $W^{1,1}(0,R_0)$, and in the sense of weak derivatives we have
\begin{align*}
	\int_{\Gamma_R \cap S} f \, |\nabla g|^{-1} \, ds = \frac{d}{dR} \, \int_{S_R} f \, da  \, .
\end{align*}
\end{lemma}
The proof is a simple application of Reynolds transport theorem and the elementary property of traces, so we shall omit it.
\begin{rem}\label{Twoformelnrem}
	Owing to Lemma \ref{Twoformeln} and to \eqref{nablagbounds}, we find for all non-negative $f \in W^{1,1}(\Omega)$ that  
	\begin{align*}
	\sqrt{2}\,\theta	\int_{\Gamma_R} f\, da \leq  \frac{d}{dR} \, \int_{C_R} f   \, dx \leq 	\frac{\sqrt{2}}{\theta}\,	\int_{\Gamma_R} f\, da   \, .
	\end{align*}
\end{rem}
Next we turn to obtaining elliptic estimates on the surfaces $\Gamma_R$.
We at first define, for $i = 1,\ldots,d$, the differential operator $\diff_i = \partial_{x_i} - \tau_i \, \tau \cdot \nabla $ acting on (weakly) differentiable functions over $\Omega$. Moreover, for the sake of simplicity in notation, we let 
\begin{align}\label{VV}
V_{\kappa} := \kappa  \tau -(\kappa \tau \cdot \tau) \,  \tau  \,\,\in\,\, \bigcap_{1\leq p < d} W^{1,p}(\Omega; \, \mathbb{R}^d) \, ,
\end{align} 
where the regularity follows from Lemma \ref{tauregu}. We also notice that $V_{\kappa}\cdot \tau = 0$ everywhere in $\Omega \setminus \{0\}$, due to the fact that $\tau$ is a unit vector. 
\begin{lemma}\label{tangentes}
	Assume that $f \in C^1(\Omega)$, $Q= 0$, and let $u \in W^{1,p}_{\text{loc}}(\Omega)$, with $\max\{p_0^{\prime},d/(d-1)\} < p < p_0$ arbitrary, be a weak solution to \eqref{P1}, \eqref{P2}. With $v_1 := \tau \cdot \nabla u \in \bigcap_{1\leq q < p^{\rm \#}} W^{1,q}_{\text{loc}}(\Omega)$ from Proposition \ref{v1prop} and with $V_{\kappa}$ from \eqref{VV}, we define on $\Omega$ an auxiliary function $$\Phi := (V_{\kappa}+\kappa\tau) \cdot \nabla v_1 + \Big[(\tau \cdot \nabla )V_{\kappa} - (V_{\kappa} \cdot \nabla)\tau + \divv (\tau) \, V_{\kappa} + \divv (\kappa \tau) \, \tau - \diff (\ln |\nabla g|)\Big]\cdot \nabla u \, .$$ Then, for all $R \in ]0,\, R_0[$, we have $\Phi\in L^2(\Gamma_R)$, $u \in W^{1,2}(\Gamma_R)$ and $v_1 \in W^{1-\frac{2}{r}}_r(S \cap \Gamma_R)$ for all $2 < r < p_0$ satisfying the integral relation (surface PDE)
\begin{align*}
\int_{\Gamma_R} \kappa \, \diff u \cdot \diff \zeta \, da +  \int_{\Gamma_R \cap S} v_1 \, [\kappa]_S\tau \cdot \nu \, \zeta\, ds = \int_{\Gamma_R} (\Phi+f) \, \zeta \, da \quad \text{ for all } \quad \zeta \in W^{1,2}(\Gamma_R) \, .
\end{align*}
\end{lemma} 
\begin{proof}
By the definition of $\diff = \nabla - \tau \, (\tau \cdot\nabla )$, we have $\nabla u = \diff u + v_1 \, \tau$ in $\Omega$. Using \eqref{VV}, with $\phi \in C^1_c(\Omega)$ arbitrary, it follows that
\begin{align}\label{tangential1}
\int_{\Omega} \kappa \nabla u \cdot \nabla \phi \, dx= \int_{\Omega} \kappa \diff u \cdot \diff \phi +V_{\kappa}\cdot \nabla u  \, \tau \cdot \nabla \phi + v_1 \, \kappa\tau \cdot \nabla \phi \, dx \, .
\end{align}
We recall that every weak solution to \eqref{P1}, \eqref{P2} belongs to $W^{2,2}_{\text{loc}}(\Omega \setminus \{0\})$. Hence, the vector fields $V_{\kappa}\cdot \nabla u  \, \tau$ and $v_1 \, \kappa\tau$ belong to $W^{1,2}_{\text{loc}}(\Omega \setminus \{0\}; \, \mathbb{R}^d)$. Restricting next to testfunctions $\phi \in C^1_c(\Omega \setminus \{0\})$, we can integrate by parts in the two last members of \eqref{tangential1}. We also note that $\tau \cdot \nu = 0$ on $S$, and the commutator formula
\begin{align*}
	\tau \cdot \nabla (V_{\kappa}\cdot \nabla u) = V_{\kappa} \cdot \nabla v_1 + [(\tau \cdot \nabla )V_{\kappa} - (V_{\kappa} \cdot \nabla)\tau]\cdot \nabla u \, .
\end{align*}
Applying these ideas we obtain that
\begin{align*}
& \int_{\Omega} (V_{\kappa}\cdot \nabla u \, \tau+ v_1 \, \kappa\tau) \cdot \nabla \phi \, dx = \int_S v_1 \, [\kappa]_S\tau \cdot \nu \, \phi \, da - \int_{\Omega} H \, \phi \, dx\\
& H:= (V_{\kappa}+\kappa\tau) \cdot \nabla v_1 + \Big[(\tau \cdot \nabla )V_{\kappa} - (V_{\kappa} \cdot \nabla)\tau + \divv (\tau) \, V_{\kappa} + \divv (\kappa \tau) \, \tau\Big]\cdot \nabla u \, .
\end{align*}
We hence obtain that $\int_{\Omega} \kappa \, \diff u \cdot \diff \phi \, dx = - \int_{S} v_1 \, [\kappa]_S\tau \cdot \nu \, \phi \, ds + \int_{\Omega} (H+f) \,  \phi \, dx$ for all $\phi \in C^1_c(\Omega \setminus \{0\})$.

Now, we choose $\phi := \zeta \, \eta \circ g$ with $\zeta\in C^2_c(\Omega)$ and $\eta \in C^1_c(0,R_0)$ arbitrary and with $g$ from \eqref{gedefin}. By the construction of $\diff$ we have $\diff g = 0$, thus
$$\int_{\Omega} \eta(g) \, \Big(\kappa \, \diff u \cdot \diff \zeta - (H+f) \,  \zeta\Big) \, dx = - \int_{S} \eta(g) \, v_1 \, [\kappa]_S\tau \cdot \nu \, \zeta \, da \, .$$ 

Due to the fact that $\eta\circ g$ is supported away from $x = 0$, it follows that $\eta(g) \, \kappa \, \diff u \cdot \diff \zeta$ and $\eta(g) \, H$ are at least of class $W^{1,1}(\Omega)$. Recall that we assume that $f \in C^1(\Omega)$. Thus, we know in particular from the proof of Proposition \ref{v1prop}, equation \eqref{ESTIGROUND}, that $w(x) = |x| \, v_1(x)$ is of class $W^{1,r}_{\text{loc}}(\Omega)$ for all $1 \leq r < p_0$. This can be used to show that $\eta(g) \, v_1 \, [\kappa]_S\tau \cdot \nu \, \zeta$ is of class $W^{1,r}(\Omega)$ for all $2 \leq r < p_0$. Hence, Lemma \ref{Twoformeln} allows to conclude that 
\begin{align*}
	\int_{0}^{R_0} \eta(R) \, \left(\int_{\Gamma_R} \frac{\kappa \, \diff u \cdot \diff \zeta - (H+f) \,  \zeta}{|\nabla g|} \, da +  \int_{\Gamma_R \cap S} \frac{v_1 \, [\kappa]_S\tau \cdot \nu\, \zeta}{|\nabla g|}\, ds \right) \, dR = 0 \, .
\end{align*}
Further, we replace $\zeta$ by $\zeta \, |\nabla g|$ and we define $\Phi := H - \diff \ln |\nabla g| \cdot \nabla u$, to obtain that
\begin{align*}
	\int_{0}^{R_0} \eta(R) \, \left(\int_{\Gamma_R} \kappa \, \diff u \cdot \diff \zeta - (\Phi+f) \,  \zeta \, da+  \int_{\Gamma_R \cap S} v_1 \, [\kappa]_S\tau \cdot \nu \, \zeta\, ds \right) \, dR = 0 \, .
\end{align*}
Since $\eta$ was an arbitrary testfunction, the claim follows.
\end{proof}
This Lemma shows that $u$ satisfies a certain elliptic transmission problem on the manifold $\Gamma_R$, with the important detail that the discontinuity of the coefficient occurs on the smooth curve $\Gamma_R \cap S$. Now, for transmission problems with interface of class $\mathscr{C}^2$, we can obtain a $W^{2,s}-$estimate on both sides for all $1< s< + \infty$ (see for instance the Theorem 6.5.2 of \cite{MR3524106} or Appendix, Lemma \ref{TransmiBasic} based on the latter).
\begin{prop}\label{curvedest}
We adopt the assumptions of Lemma \ref{tangentes}. Then, for all $1 < s < +\infty$,
\begin{align*}
	\|D^2 u\|_{L^s(\Gamma_R \setminus S)} \leq C \, (\|f\|_{L^s(\Gamma_{R})} + \|\nabla v_1\|_{L^s(\Gamma_{R})} + R^{-1} \, \|\nabla u\|_{L^s(\Gamma_{R})} + R^{-2} \, \|u\|_{L^s(\Gamma_{R})}) \, ,
\end{align*}
with the constant being independent of $u$ and $0< R \leq R_0$.
\end{prop}
\begin{proof}
	Choosing $0 < R \leq R_0/\sqrt{2}$ arbitrary, we at first locally flatten the manifold $\Gamma_R$. To this aim we recall the definition \eqref{phinull} and we fix two numbers $b_0 < b_1 $ such that
	\begin{align}\label{bnot}
		\frac{1}{\sqrt{1+\theta \, (\phi_0(R_0))^2}} < b_0 < b_1 < 1 \, .
	\end{align}
Let $M$ be the function introduced in Lemma \ref{Mandmu}. 
In the upper half-space $\{(\bar{x},x_d) \, : \, x_d >0\}$, the surface $\Gamma_R$ is described as the graph of the function $\gamma: \, B_{\sqrt{2}R}^{d-1}(0) \rightarrow [0, \, M^{-1}_+(R^2)]$ given as
\begin{align}\label{gammadef}
	\gamma(\bar{x}) := M^{-1}_+(R^2 - |\bar{x}|^2/2) \quad \text{ for } \quad |\bar{x}| \leq \sqrt{2}R \, .
\end{align}
Hereby, $B^{d-1}_{\rho}(0)$ denotes the $(d-1)$-dimensional unit ball with radius $\rho$ and centre $0$. Since $M^{-1}_+$ is increasing, the bounds proved in Lemma \ref{Mandmu} show that, when restricting $\gamma$ to $B_{\sqrt{2}b_1R}^{d-1}(0)$, we have
\begin{align}\label{goulde}
& \sqrt{2(1-b_1^2)}\,  R \leq	M^{-1}_+((1-b_1^2)\,R^2) \leq \gamma(\bar{x}) \leq M^{-1}_+(R^2) \leq \sqrt{\frac{2}{\theta}} \, R \\
& \qquad \text{ for } \quad |\bar{x}| \leq b_1 \, \sqrt{2} \,  R \, .\nonumber
\end{align}
Using \eqref{gammadef}, we can further compute that $\partial_{x_i} \gamma(\bar{x}) = \frac{-x_i}{M^{\prime}(\gamma(\bar{x}))}$ for $i = 1,\ldots,d-1$. Recalling \eqref{mulineargrowth}, we have $M^{\prime}(\gamma(\bar{x}))/\gamma(\bar{x}) \geq \theta$, and the inequalities \eqref{goulde} yield
\begin{align}\label{gradientbounded}
	|\nabla_{\bar{x}} \gamma(\bar{x})| \leq \frac{b_1}{\theta\, \sqrt{1-b_1^2}}  \quad \text{ for } \quad |\bar{x}| \leq b_1 \, \sqrt{2} \,  R \, .
\end{align}
Next we compute that
\begin{align*}
	\partial^2_{x_i,x_j} \gamma(\bar{x}) = -\frac{1}{M^{\prime}(\gamma)} \, \Big(\delta_{ij}+ \frac{x_i \, x_j}{(M^{\prime}(\gamma))^2} \, M^{\prime\prime}\Big) \quad \text{ for } \quad i,j = 1,\ldots,d-1\, ,
\end{align*}
and by similar means as for the gradient, we see that
\begin{align}\label{hessianbounded}
	|D^2_{\bar{x}} \gamma(\bar{x})|_{\infty} \leq \Big(1+\frac{b_1^2}{\theta^2(1-b_1^2)}\Big) \, \frac{1}{\theta \,\sqrt{2(1-b_1^2)} \, R}  \quad \text{ for } \quad |\bar{x}| \leq b_1 \, \sqrt{2} \,  R \, .
\end{align}
We note that by the result of Lemma \ref{curve} and the choice \eqref{bnot} of $b_0$ we find for $0 < R \leq R_0/\sqrt{2}$ that the radius $\hat{\rho}(R)$ of the intersection $\Gamma_R \cap S$ satisfies
\begin{align}\label{inclusion}
	\hat{\rho}(R) \leq \frac{\sqrt{2} \, R}{(1+\theta\, \phi_0^2(\sqrt{2}R))^{\frac{1}{2}}} < b_0 \, \sqrt{2}\, R \, \quad \Rightarrow \quad B_{\hat{\rho}(R) }^{d-1}(0)  \subset\!\!\subset B_{b_0\sqrt{2}R}^{d-1}(0) \, ,
\end{align}
Next for $\bar{x} \in B_{\sqrt{2}R}^{d-1}(0)$ and for a function $f$ defined in $\Omega$, we define $\tilde{f}(\bar{x})	:= f(\bar{x}, \,\gamma(\bar{x}))$ to be the constant extension in normal direction. For a vector $a = (a_1, \ldots ,a_d)$ we also recall the notation $\bar{a} = (a_1,\ldots,a_{d-1})$, and we then introduce 
\begin{align*}
\mathscr{D}(\bar{x}) := \sqrt{1+|\nabla \gamma(\bar{x})|^2}, \quad \text{and} \quad Q(\bar{x}) := [e^1, \ldots, e^{d-1}]
- \tilde{\tau}(\bar{x}) \otimes \tilde{\bar{\tau}}(\bar{x}) \in \mathbb{R}^{d\times (d-1)} \, ,
\end{align*}
where $[e^1, \ldots, e^{d-1}]$ is the matrix with columns consisting of the $d-1$ first standard basis vectors of $\mathbb{R}^d$. Recalling the characterisation of the curve $\Gamma_R \cap S$ in Lemma \ref{curve}, the result of Lemma \ref{tangentes} now implies that
\begin{align}\label{transmi2d}
	& \int_{B^{d-1}_{b_1 \sqrt{2}R}} \mathscr{D} \, Q^{\sf T}\, \kappa\, Q \, \nabla_{\bar{x}} \tilde{u} \cdot \nabla_{\bar{x}} \tilde{\zeta} \, d\bar{x} \nonumber\\ 
	&=  -\int_{\partial B^{d-1}_{\hat{\rho}(R)}} \widetilde{v_1} \, [\kappa]\tilde{\tau}\cdot \widetilde{\nu} \, \tilde{\zeta} \, ds 
	 +	\int_{B^{d-1}_{b_1 \sqrt{2}R}} \mathscr{D} \, (\tilde{\Phi}+\tilde{f}) \,  \tilde{\zeta} \, d\bar{x} \quad \text{ for all } \quad \tilde{\zeta} \in C^1_c(B^{d-1}_{b_1 \sqrt{2}R}) \, .
\end{align}	
We introduce the abbreviation
\begin{align*}
A(\bar{x}) :=  \mathscr{D}(\bar{x}) \, Q^{\sf T}(\bar{x})\,\kappa Q(\bar{x})\, \quad \text{ for } \quad \bar{x} \in B^{d-1}_{b_1 \sqrt{2}R}\setminus \partial B^{d-1}_{\hat{\rho}(R)} \, ,
\end{align*}
and see that the map $\bar{x} \mapsto A(\bar{x})$ is continuously differentiable in $B^{d-1}_{b_1 \sqrt{2}R}\setminus \partial B^{d-1}_{\hat{\rho}(R)}$. By means of \eqref{gradientbounded} and \eqref{hessianbounded}, we moreover easily verify that
\begin{align}\label{properA}
\sup_{\bar{x} \in B^{d-1}_{b_1 \sqrt{2}R} \setminus \partial B^{d-1}_{\hat{\rho}(R)}} |A(\bar{x})|_{\infty} \leq \frac{C(b_1)}{\theta}, \quad \sup_{\bar{x} \in B^{d-1}_{b_1 \sqrt{2}R} \setminus \partial B^{d-1}_{\hat{\rho}(R)}} |\nabla_{\bar{x}}A(\bar{x})|_{\infty}\leq  \frac{C(b_1)}{\theta^3 \, R}\, .
\end{align}
Moreover, invoking that $\tau = \nabla g/|\nabla g|$ is normal on $\Gamma_R$, we have $\tilde{\tau}_d(\bar{x}) = 1/\mathscr{D}(\bar{x})$, and we can show for $\eta\in \mathbb{R}^{d-1}$ arbitrary that
\begin{align*}
	|Q(\bar{x}) \, \eta|^2 = |\eta|^2 - (\tilde{\bar{\tau}}(\bar{x}) \cdot \eta)^2 \geq \Big(\frac{1}{\mathscr{D}(\bar{x})}\Big)^2 \,  |\eta|^2\, .
\end{align*}
Hence, invoking \eqref{gradientbounded} once again,
\begin{align*}
A(\bar{x}) \eta \cdot \eta = & \mathscr{D}(\bar{x}) \, \kappa \,  Q(\bar{x})\eta \cdot Q(\bar{x})\eta \geq  \mathscr{D}(\bar{x}) \,k_0 \, |Q(\bar{x}) \, \eta|^2\\
\geq & \frac{\theta\, k_0}{\sqrt{1+\theta^2}} \, |\eta|^2 \quad \text{ for all } \quad \eta \in \mathbb{R}^2, \, \bar{x} \in B^{d-1}_{b_1 \sqrt{2}R} \setminus \partial B^{d-1}_{\hat{\rho}(R)} \, .
\end{align*}
Hence, \eqref{transmi2d} is the weak form of a $(d-1)$-dimensional transmission problem with coefficient matrix $A$ and interface $ \partial B^{d-1}_{\hat{\rho}(R)}$.

Recall that we for simplicity assumed $f \in C^{1}(\Omega)$. For $R > 0$, the distance of $\Gamma_{R}$ to zero is strictly positive, hence the function $\Phi$ as well as the vectors $\tau$ and $\nu$ are smooth therein. The function $\tilde{v}_1$ belongs to $W^{1,r}(\Omega \setminus B_{\epsilon}(0))$ for all $0 < \epsilon \ll R$. Thus, via the trace theorem, and since the function $\gamma$ is of class $C^{1,1}$, we find that $\tilde{v}_1$ belongs to $W^{1-1/r}_r(B^{d-1}_{b_1 \sqrt{2}R})$ and to $W^{1-2/r}_r(\partial B^{d-1}_{\hat{\rho}(R)})$ for all $2 < r < +\infty$. Hence, via the Sobolev embedding theorem, we can choose $r \geq 2d/(d-1)$ to see that $\tilde{v}_1 \in W^{1-2/r}_r(\partial B^{d-1}_{\hat{\rho}(R)}) \subset W^{1/2}_2(\partial B^{d-1}_{\hat{\rho}(R)})$.
Due to this regularity, of the right-hand side and the surface term, we directly obtain that $\tilde{u}$ satisfying \eqref{transmi2d} belongs to $W^{2,2}_{\text{loc}}(B^{d-1}_{b_1 \sqrt{2}R} \setminus \partial B^{d-1}_{\hat{\rho}(R)} )$ and is thus a strong solution.

In order to remove the constraint on testfunctions in \eqref{transmi2d}, we fix  $\zeta_0 \in C^2_c(B^{d-1}_{b_1 \sqrt{2}R})$ which is a non-negative function bounded by one, equal to one on $B_{b_0 \, \sqrt{2}R}^{d-1}$, and the derivatives of which satisfy $|\nabla_{\bar{x}} \zeta_0| \leq C_0\, (b_1-b_0)^{-1} \, R^{-1}$ and $|D^2_{\bar{x}}\zeta_0| \leq C_1 \, (b_1-b_0)^{-2} \, R^{-2}$. Then with $w := \tilde{u}\zeta_0$ and $\psi = \widetilde{v_1} \zeta_0$ we get after a few calculations
\begin{align}\label{transmi2dbis}
	\int_{B^{d-1}_{b_1 \sqrt{2}R}} A(\bar{x}) \, \nabla_{\bar{x}} w \cdot \nabla_{\bar{x}} \zeta \, d\bar{x} = -\int_{\partial B^{d-1}_{\hat{\rho}(R)}} \psi \, [\kappa]\tilde{\tau}\cdot \tilde{\nu}\, \zeta \, ds + \int_{B^{d-1}_{b_1 \sqrt{2}R}} \Psi(\bar{x}) \,  \zeta \, d\bar{x}  \\
\text{ for all } \quad \zeta \in C^1(\overline{B^{d-1}_{b_1 \sqrt{2}R}}) \quad \text{ with } \quad \Psi =\mathscr{D} \, (\tilde{\Phi}+\tilde{f}) - 2 \ A \,\nabla_{\bar{x}}\tilde{u} \cdot \nabla_{\bar{x}}\zeta_0 - \divv_{\bar{x}}(A \, \nabla_{\bar{x}} \zeta_0) \, \tilde{u}\,. \nonumber
\end{align}	
The next step is to obtain weaker estimates on the second derivatives. In order to track dependence on $R$, we first re-scale the problem by introducing $y = \bar{x}/\hat{\rho}(R)$, in such a way that the interface $\partial B^{d-1}_{\hat{\rho}(R)}$ is mapped on a sphere of unit length. A useful abbreviation is $k(R):= b_1\sqrt{2}R/\hat{\rho}(R)$. Owing to Lemma \ref{curve}, and the choices \eqref{bnot} of the parameters $b_0 < b_1$, we see that $k(R) \geq b_1/b_0 > 1$ for $R \leq R_0/\sqrt{2}$. In fact, in case of the condition (A2), we have $k(R) \geq b_1 \, (1+\theta\, \phi_0^2(\sqrt{2}R))^{1/2}$, and $k(R)$ tends to infinity as $R \rightarrow 0$, while it remains bounded from below under the condition (A2').

A function $f$ defined in $B^{d-1}_{b_1 \sqrt{2}R}$ transforms to a function $\hat{f}(y) := f(\hat{\rho}(R) \, y)$ defined in $B^{d-1}_{k(R)}(0)$, and we have $\nabla_{\bar{x}} = \hat{\rho}(R)^{-1} \, \nabla_y$. Applying the transformation formula to \eqref{transmi2dbis}, we get
\begin{align*}
	\int_{B^{d-1}_{k(R)}} \hat{A}(y)\, \nabla_y \hat{w} \cdot \nabla_y \hat{\zeta} \, dy = & - \hat{\rho}(R) \, \int_{\partial B^{d-1}_{1}(0)} \hat{\psi} \, \widehat{[\kappa]\tilde{\tau}\cdot \widetilde{\nu}}\, \hat{\zeta} \, ds \\
	&+ (\hat{\rho}(R))^2 \,	\int_{B^{d-1}_{k(R)}} \hat{\Psi}(y) \,  \hat{\zeta} \, dy \quad \text{ for all } \quad \hat{\zeta} \in C^1(\overline{B^{d-1}_{k(R)}}) \, .
\end{align*}	
As shown in Appendix, Lemma \ref{TransmiBasic}, this transmission problem allows, for all $1 < s < +\infty$, for a bound
\begin{align*}
	\|D^2_y \hat{w}\|_{L^s(B^{d-1}_{k(R)})} \leq  c(s) \, \Big( & \|\nabla_y \hat{A}\|_{L^{\infty}} \, \|\nabla_y\hat{w}\|_{L^s(B^{d-1}_{k(R)})} + \hat{\rho}(R) \, \|\widehat{[\kappa]\tilde{\tau}\cdot \widetilde{\nu}}\, \hat{\psi}\|_{W^{1-\frac{1}{s}}_s(\partial B^{d-1}_1)} \\
	& + (\hat{\rho}(R))^2 \,  \|\hat{\Psi}\|_{L^s(B^{d-1}_{k(R)})}\Big) \, .
\end{align*}
Moreover, \eqref{properA} and the transformation rules allow to estimate $\|\nabla_y \bar{A}\|_{L^{\infty}} \leq C \,  \hat{\rho}(R)/ (\theta^{3} \,R)$, hence
\begin{align*}
	\|D^2_y \hat{w}\|_{L^s(B^{d-1}_{k(R)})} \leq 
	 c(s) \, \Big(& \frac{C}{\theta^3 \, k(R)} \, \|\nabla_y\hat{w}\|_{L^s(B^{d-1}_{k(R)})}  + \hat{\rho}(R) \, \|\widehat{[\kappa]\tilde{\tau}\cdot \widetilde{\nu}}\, \hat{\psi}\|_{W^{1-\frac{1}{s}}_s(\partial B^{d-1}_1)}\\
	 & + (\hat{\rho}(R))^2 \,  \|\hat{\Psi}\|_{L^s(B^{d-1}_{k(R)})}\Big) \, . 
\end{align*}
Observe that $\tilde{\tau}$ and $\tilde{\nu}$ are orthogonal on $\partial B^{d-1}_{\hat{\rho}(R)}$ and thus $[\kappa] \tilde{\tau}\cdot \tilde{\nu} =  [\kappa^\circ] \tilde{\tau}\cdot \tilde{\nu}$ with the trace-less part $\kappa^\circ$ of $\kappa$. As a consequence of Corollary \ref{w1infty}, we then have
\begin{align*}
	|[\kappa]\tilde{\tau}(\bar{x})\cdot \widetilde{\nu}(\bar{x}) - [\kappa]\tilde{\tau}(\xi^{\prime})\cdot \widetilde{\nu}(\xi^{\prime})| \leq 2\lambda_1\, |[\kappa^0]|_2 \, \frac{|\bar{x}-\xi^{\prime}|}{\hat{\rho}(R)} \, \text{ for  all } \bar{x}, \, \xi^{\prime} \in \partial B^{d-1}_{\hat{\rho}(R)} \, ,
\end{align*}
which implies that
\begin{align}\label{solle}
	|\widehat{[\kappa]\tilde{\tau}\cdot \tilde{\nu}}(y) - \widehat{[\kappa]\tilde{\tau}\cdot \tilde{\nu}}(y^*)| \leq 2\lambda_1\, |[\kappa^\circ]|_2 \, |y-y^*| \, \text{ for  all } y, \, y^* \in \partial B_{1}^{d-1} \, .
\end{align}
Hence, we conclude that
\begin{align*}
\|\widehat{[\kappa]\tilde{\tau}\cdot \tilde{\nu}}\, \hat{\psi}\|_{W^{1-\frac{1}{s}}_s(\partial B^{d-1}_1)} \leq & \|\overline{[\kappa]\tilde{\tau}\cdot \tilde{\nu}}\|_{L^{\infty}(\partial B^{d-1}_1)}\, \|\hat{\psi}\|_{W^{1-\frac{1}{s}}_s(\partial B^{d-1}_1)} + 4\pi \, \lambda_1\, |[\kappa^\circ]|_2 \, \|\hat{\psi}\|_{L^s(\partial B^{d-1}_1)} \\
\leq & C \,  |[\kappa^\circ]|_2 \, \|\hat{\psi}\|_{W^{1-\frac{1}{s}}_s(\partial B^{d-1}_1)} \, .
\end{align*}
Employing further the fact that $W^{1-1/s}_s(\partial B^{d-1}_1)$ can be identified as the trace space of $W^{1,s}(\mathbb{R}^{d-1})$ on $\partial B^{d-1}_1$, with trace embedding constant $c_0$, we get $\|\widehat{[\kappa]\tilde{\tau}\cdot \tilde{\nu}}\, \hat{\psi}\|_{W^{1-1/s}_s(\partial B^{d-1}_1)} \leq  C_1 \,  |[\kappa^\circ]|_2 \, \|\nabla_y \hat{\psi}\|_{L^{s}(\mathbb{R}^{d-1})}$. Hence, since $\hat{\psi}$ as compact support in $B^{d-1}_{k(R)}$, it follows that
\begin{align*}
	\|D^2_y \hat{w}\|_{L^s(B^{d-1}_{k(R)})} \leq 
	c(s) \, \Big( & \frac{C}{\theta^3 \, k(R)} \, \|\nabla_y\hat{w}\|_{L^s(B^{d-1}_{k(R)})}  + C_1 \,  |[\kappa^\circ]|_2 \, \hat{\rho}(R) \, \|\nabla_y \hat{\psi}\|_{L^{s}(B^{d-1}_{k(R)})}\\
	& + (\hat{\rho}(R))^2 \,  \|\hat{\Psi}\|_{L^s(B^{d-1}_{k(R)})}\Big) \, . 
\end{align*}
We transform back from $y$ to $\bar{x} = \hat{\rho}(R) \, y$ and get
\begin{align*}
\|D^2_{\bar{x} } w\|_{L^s(B^{d-1}_{b_1 \sqrt{2}R})} \leq  c_s \, \Big(& \frac{C}{\theta^3\,R} \, \|\nabla_{\bar{x}} w\|_{L^s(B^{d-1}_{b_1 \sqrt{2}R})} + C_1 \,  |[\kappa^{\circ}]|_2 \,  \|\nabla_{\bar{x}} \psi\|_{L^{s}(B^{d-1}_{b_1 \sqrt{2}R})}+\|\Psi\|_{L^s(B^{d-1}_{b_1 \sqrt{2}R})}\Big) \, . 
\end{align*}
We might rescale $R \approx b_0 \sqrt{2}R$. After a few elementary steps, where we use the bounds for the derivatives of $\zeta_0$, we arrive at the following form of the estimate, valid for all $0 < R \leq b_0R_0$ with $\chi = b_1/b_0 > 1$:
\begin{align*}
	\|D^2_{\bar{x}}\tilde{u}\|_{L^s(B^{d-1}_{R})} \leq C \, \Big(\|\tilde{f}\|_{L^s(B^{d-1}_{\chi R})} + \|\widetilde{\nabla v_1}\|_{L^s(B^{d-1}_{\chi R})} + \frac{1}{R} \, \|\widetilde{\nabla u}\|_{L^s(B^{d-1}_{\chi R})} + \frac{1}{R^2} \, \|\tilde{u}\|_{L^s(B^{d-1}_{\chi R})}\Big) \, .
\end{align*}
It is readily seen that going back to curved coordinates does not change the structure of the estimate. With $\Gamma_R^{\prime} := \{(\bar{x}, \, \gamma(\bar{x}))\, ; \, |\bar{x}| < b_0 \, \sqrt{2}\,R\}$, we hence obtain that
\begin{align*}
	\|\diff^2 u\|_{L^s(\Gamma_R^{\prime})} \leq C \, \Big(\|f\|_{L^s(\Gamma_{R})} + \|\nabla v_1\|_{L^s(\Gamma_{R})} + \frac{1}{R} \, \|\nabla u\|_{L^s(\Gamma_{R})} + \frac{1}{R^2} \, \|u\|_{L^s(\Gamma_{R})}\Big) \, .
\end{align*}
Note that the surface $\Gamma_{R} \setminus \Gamma_R^{\prime}$ does not intersect the surface $S$, owing to \eqref{inclusion}.
The integral relation of Lemma \ref{tangentes} with testfunctions $\zeta \in C_c(\Gamma_R \setminus \overline{\Gamma}_{R}^{\prime})$ is therefore a transmission problem with coefficient of class $C^1$.
After another similar localising argument, we therefore obtain a global bound on $\Gamma_R$ of the form
\begin{align*}
	\|\diff^2 u\|_{L^s(\Gamma_R)} \leq C \, \Big(\|f\|_{L^s(\Gamma_{R})} + \|\nabla v_1\|_{L^s(\Gamma_{R})} + \frac{1}{R} \, \|\nabla u\|_{L^s(\Gamma_{R})} + \frac{1}{R^2} \, \|u\|_{L^s(\Gamma_{R})}\Big) \, .
\end{align*}
Finally, we can obtain the full Hessian using the identity $\nabla u = \diff u + v_1 \, \tau$. After a few calculations, we see that $|D^2 u| \leq |\diff^2 u| + |\nabla v_1| + 2\, |\nabla \tau| \, |\nabla u|$, and since $|\nabla \tau| \leq c \, R^{-1}$ on $\Gamma_R$, the claim follows.
\end{proof}
\begin{lemma}\label{SecondDeriv}
	Let $\max\{p_0^{\prime},d/(d-1)\} < r < p_0^*$ and $1 \leq q < r^{\rm \#}$. Then, for every $s_1 > p_0^{\prime} $ and $\beta_1$ as in Proposition \ref{v1prop}, we find a constant $c$ such that
	\begin{align*}
		\|D^2u\|_{L^q(\Omega^{\prime\prime} \setminus S)} \leq c \, (\|f\|_{L^q(\Omega^{\prime})} + \|E_{\beta_1} f\|_{L^{s_1}(\Omega^{\prime})}  +\|u\|_{W^{1,r}(\Omega^{\prime})}) \, .
	\end{align*}
\end{lemma}
\begin{proof}
	We raise the inequality of Proposition \ref{curvedest} to the power $s$. Recalling the remark \ref{Twoformelnrem}, we obtain for $R \leq R_0$ arbitrary that
	\begin{align}\begin{split}\label{roulez}
		\frac{d}{dR} \int_{C_R} |D^2u|^s \, dx \leq \frac{C}{\theta} \, \Big( &	\frac{d}{dR} \int_{C_R} |f|^s + |\nabla v_1|^s \, dx \\
		&  + \frac{1}{R^s} \,\frac{d}{dR} \int_{C_R} |\nabla u|^s \, dx  + \frac{1}{R^{2s}} \,\frac{d}{dR} \int_{C_R} |u|^s \, dx  \Big)\, .\end{split}
	\end{align}
	We want to integrate this inequality. We let $0 < \epsilon < R_0$, and we notice that 
	\begin{align*}
		\int_{\epsilon}^{R_0} \frac{1}{R^s} \,\frac{d}{dR} \int_{C_R} |\nabla u|^s \, dxdR = & \frac{1}{R_0^s} \,\int_{C_{R_0}} |\nabla u|^s \, dx - \frac{1}{\epsilon^s} \,\int_{C_{\epsilon}} |\nabla u|^s \, dx \\
		& + s \, 	\int_{\epsilon}^{R_0} \frac{1}{R^{1+s}} \, \int_{C_R} |\nabla u|^s \, dxdR \, .
	\end{align*}
	We let $1 < s < r$. By means of H\"older's inequality, $ \int_{C_R} |\nabla u|^s \, dx \leq \|\nabla u\|_{L^{r}(C_{R_0})}^s \, |C_R|^{1-\frac{s}{r}}$. We can easily verify that $C_R \subset B_{\sqrt{2\theta^{-1}}R}(0)$, hence $|C_R| \leq 2^{\frac{d}{2}} \, \theta^{-\frac{d}{2}} \, \omega_d \, R^d$. For all $1 < s < dr/(d+r)$, it hence follows that
	\begin{align*}
		\frac{1}{R^{1+s}}	\, |C_R|^{1-\frac{s}{r}} \leq c_1 \, R^{d\Big[1-\frac{s}{r}\Big] - s - 1} = c \, R^{\delta_0-1} \quad \text{ with } \quad \delta_0 = d\Big[1-\frac{s}{r}\Big]- s > 0\, , 
	\end{align*}
	with $c_1 := 2^{\frac{d}{2}} \, \theta^{-\frac{d}{2}} \, \omega_d$. Hence
	\begin{align*}
		&	\int_{\epsilon}^{R_0} \frac{1}{R^{1+s}} \, \int_{C_R} |\nabla u|^s \, dxdR \leq \frac{c_1}{\delta_0} \, \|\nabla u\|_{L^{r}(C_{R_0})}^s \, 	(R_0^{\delta_0}-\epsilon^{\delta_0}) \, ,\\
		\text{ and }\quad &  \frac{1}{\epsilon^s} \,\int_{C_{\epsilon}} |\nabla u|^s \, dx\leq c_1 \, \epsilon^{\delta_0} \|\nabla u\|_{L^{r}(C_{R_0})}^s\, .
	\end{align*}
	This allows to show that
	\begin{align*}
		\limsup_{\epsilon \rightarrow 0} \int_{\epsilon}^{R_0} \frac{1}{R^s} \,\frac{d}{dR} \int_{C_R} |\nabla u|^s \, dxdR \leq\frac{1}{R_0^s} \,\int_{C_{R_0}} |\nabla u|^s \, dx  + \frac{c_1 \, s\, R_0^{\delta_0}}{\delta_0} \,  \|\nabla u\|_{L^{r}(C_{R_0})}^s\, .
	\end{align*}
	For $\epsilon > 0$, we next notice in the same spirit that 
	\begin{align*}
		\int_{\epsilon}^{R_0} \frac{1}{R^{2s}} \,\frac{d}{dR} \int_{C_R} |u|^s \, dxdR = & \frac{1}{R_0^{2s}} \,\int_{C_{R_0}} |u|^s \, dx - \frac{1}{\epsilon^{2s}} \,\int_{C_{\epsilon}} |u|^s \, dx\\
		& + 2s \, 	\int_{\epsilon}^{R_0} \frac{1}{R^{1+2s}} \, \int_{C_R} |u|^s \, dxdR \, .
	\end{align*}
	Now $ \int_{C_R} |u|^s \, dx \leq \|u\|_{L^{r^*}(\Omega)}^s \, |C_R|^{1-\frac{s}{r^*}}$ with $r^* = dr/(d-r)$. Since $s < dr/(d+r)$, we also have
	\begin{align*}
		\frac{1}{R^{1+2s}}	\, |C_R|^{1-\frac{s}{r^*}} \leq c_1 \, R^{d\Big[1-\frac{s}{r^*}\Big] - 1 - 2s} = c_1 \, R^{\delta_0-1} \, . 
	\end{align*}
We choose $s := q < dr/(d+r)$ and integrate \eqref{roulez} over $(0,R_0)$ to obtain that
\begin{align}
		\int_{C_{R_0}} |D^2u|^q \, dx \leq c \, \Big( & \int_{C_{R_0}} |f|^q + |\nabla v_1|^q \, dx  +\int_{C_{R_0}} |\nabla u|^r + |u|^r \, dx  \Big)\, .
\end{align}
We invoke Proposition \ref{v1prop} to estimate $\|\nabla v_1\|_{L^q}$, and the claim follows.
\end{proof}
The following statement completes the proof of Theorem \ref{main}.
\begin{coro}\label{SecondDerivcoro}
	Let $\max\{p_0^{\prime},d/(d-1)\} < p < p_0$ and $1 \leq q < dp/(d+p)$. Then
	\begin{align*}
		\|D^2u\|_{L^q(\Omega^{\prime\prime} \setminus S)} \leq c \, (\|f\|_{L^{p^{\rm \#}}(\Omega^{\prime})} + \|E_1 f\|_{L^p(\Omega^{\prime})}  +\|u\|_{L^{p}(\Omega^{\prime})}) \, .
	\end{align*}
	\end{coro}
\begin{proof}
In Lemma \ref{SecondDeriv}, we choose $r = p$, $s_1 = p$ and $\beta_1 = 1$. We invoke \eqref{ugroeg} to estimate $\|u\|_{W^{1,p}}$, and the claim follows.
\end{proof}
In the next section, we shall discuss the higher regularity result of Theorem \ref{main2}. In order to introduce to this topic, let us exhibit an elementary consequence of Lemma \ref{SecondDeriv}.
\begin{rem}\label{trivialcase}
Suppose that the constant $p_0$ of condition \eqref{Meyers} satisfies $p_0  = + \infty$ and the assumptions of Lemma \ref{SecondDeriv} are true with $s_1 = + \infty$. Then we can choose arbitrary $p < + \infty$ in the latter statement and obtain the integrability of second derivatives for all $q < p^{\rm \#} = d$. 
\end{rem}

\section{Improving regularity}\label{gnsection}

In this section, we discuss the proof of Theorem \ref{main2}.
Relying on the H\"older continuity of solutions to \eqref{P1}, \eqref{P2} (which follows from the property \eqref{deGiorgiNash}) and on the regularity of second derivatives obtained in Lemma \ref{SecondDeriv}, a Gagliardo-Nirenberg interpolation inequality yields the higher integrability of the gradient. In the case of (A2'), both domains $\Omega_1$ and $\Omega_2$ separated by $S$ are Lipschitzian. For every domain $G = \Omega^{\prime} \cap \Omega_i$, where $\Omega^{\prime} \subset\!\!\subset \Omega$ is an arbitrary smooth subset, we can rely on the following inequality (see \cite{nirenberg66b}, Theorem 1', Appendix, Lemma \ref{scaleGN}):
\begin{align}\label{GN}
\|\nabla u\|_{L^{q/a}} \leq C_1 \, \|D^2 u\|_{L^{q}}^a \, [u]_{C^{\lambda}}^{1-a} + C_2 \, [u]_{C^{\lambda}} \, ,
\end{align}
where all norms are taken over the domain $G$ and the constants $C_1$, $C_2$ depend only on $G$, while the parameters are only subject to the following restrictions:
\begin{align}\label{RESTRI1}
0 \leq \lambda \leq 1, \quad 1 \leq q < +\infty, \quad \begin{cases} a \in [\frac{1-\lambda}{2-\lambda}, \, 1] & \text{ for } q \neq d\, ,\\
	a \in [\frac{1-\lambda}{2-\lambda}, \, 1[ &  \text{ for } q = d	\, .
	\end{cases}
\end{align}
Using this inequality and the property \eqref{deGiorgiNash}, we can improve the regularity iteratively. Recall that \eqref{deGiorgiNash} guarantees that there is $0 < \alpha_0 \leq  1$ such that every local weak solution to \eqref{P1}, \eqref{P2} with $f = 0$ and $Q = 0$ satisfies $u \in C^{\alpha_0}_{\text{loc}}(\Omega)$. Moreover, if $f \neq 0$, and $f \in L^{s_0}$ for $s_0> d/2$, we obtain from the condition \eqref{deGiorgiNash} the H\"older regularity $u \in C^{\lambda_0}_{\text{loc}}(\Omega)$ where
\begin{align}\label{lambda0}
	\lambda_0 := \min\Big\{\alpha_0, \, 1-\frac{d}{s_0^*}\Big\} \, .
\end{align}
 If $u$ is a local weak solution to \eqref{P1}, \eqref{P2} it then obeys
\begin{align}\label{unash}
	\|u\|_{C^{\lambda_0}(\Omega^{\prime\prime})} \leq c(k_0,k_1,\Omega^{\prime\prime}) \, (\|f\|_{L^{s_0}(\Omega^{\prime})}+\|u\|_{L^{s_0^*}(\Omega^{\prime})}) \quad \text{ for all } \quad \Omega^{\prime\prime} \subset\!\!\subset \Omega^{\prime} \subset\!\!\subset  \Omega \,.
\end{align}
\begin{lemma}\label{Lippeschitz}
We assume that $\Omega$ satisfies (A1), (A2') and (A3) and that $u$ is a local weak solution to \eqref{P1}, \eqref{P2}. Let $f \in L^{s_0}(\Omega)$ with $s_0 > d/2$, and assume that there are $s_1 > p_0^{\prime}$, $s_1 \geq s_0$ and $0 \leq \beta_1 \leq 1$ such that $\|E_{\beta_1} f\|_{L^{s_1}(\Omega)} < + \infty$. With $m_0 = \min\{p_0,s_1\}$, which we assume satisfying $m_0 < +\infty$ (cf.\ Rem.\ \ref{trivialcase}, for $m_0 = + \infty$), we define
\begin{align*}
	\bar{r} := d \, \min\Big\{\Big(\frac{m_0}{d+(\beta_1-1)\, m_0}\Big)^\circ, \, \Big(\frac{1}{1-\lambda_0}\Big)^\circ\Big\}\, ,
\end{align*}
where for $x > 0$ and $y \in \mathbb{R}$, we set $(x/y)^{\circ} := x/y$ for $y >0$, $(x/y)^\circ$ arb.\ large positive for $y =0$, $(x/y)^{\circ} := +\infty$ for $y < 0$, and where $0 < \lambda_0 \leq 1$ is the constant of condition \eqref{lambda0}.
Then $\nabla u \in L^r_{\text{loc}}(\Omega; \, \mathbb{R}^d)$ for all $1 \leq r < \bar{r}$ if $\bar{r}$ is finite and for $r = +\infty$ otherwise. Moreover, $D^2 u \in L^q_{\text{loc}}(\Omega\setminus S; \, \mathbb{R}^{d\times d})$ for all $1 \leq q < \bar{r}d/(\bar{r}+d)$. For all $\Omega^{\prime\prime}\subset\!\!\subset \Omega^{\prime}\subset\!\!\subset\Omega$, there is a constant $c = c(r,q,\Omega^{\prime\prime},k_0,k_1,s_0,s_1,\beta_1)$ such that
\begin{align*}
	\|D^2 u\|_{L^{q}(\Omega^{\prime\prime}\setminus S)} + 	\|\nabla u \|_{L^{r}(\Omega^{\prime\prime})} \leq c \, (\|E_{\beta_1} f\|_{L^{s_1}(\Omega^{\prime})} +\|f\|_{L^{s_0}(\Omega^{\prime})} + \|u\|_{L^{s_0^*}(\Omega^{\prime})})\, .
	\end{align*}
\end{lemma}
Note: In particular, \eqref{lambda0} implies that $\bar{r} \leq s_0^*$ and $d\bar{r}/(d+\bar{r}) \leq s_0$, so that for $q < d\bar{r}/(d+\bar{r}) $, we have $\|f\|_{L^q} \leq c \, \|f\|_{L^{s_0}}$.
\begin{proof}
If $\bar{r}$ is infinite, then $\lambda_0 = 1$. Hence $\nabla u \in L^{\infty}$, and the claim for the second derivatives follows from Lemma \ref{SecondDeriv} and the preliminary note. Thus, we assume that $\bar{r}$ is finite. 
	
We let $a_0 := d/(d+\bar{r})$ and we fix a $\max\{p_0^{\prime},d/(d-1)\} < r_0<\min\{p_0,\, \bar{r}\}$. These choices imply that $(r_0/d+1) \, \max\{a_0, \, 1/r_0\} < 1$. We choose any
\begin{align*}
	\max\Big\{\frac{r_0+d}{dr_0}, \, a_0\Big(1+\frac{r_0}{d}\Big) \Big\} < \chi < 1 \, ,
\end{align*}
and we define $q_0 := \chi \, r_0d/(r_0+d)$ and 
\begin{align*}
r_m := \frac{1}{a_0} \, q_{m-1} \quad \text{ and } \quad q_{m} = \chi \, \frac{r_md}{r_m+d}	\quad \text{ for } m \in \mathbb{N} \, .
\end{align*}
We can easily verify that
\begin{align*}
	r_m = \frac{\chi \, d}{a_0 \, (r_{m-1} + d)}  \, r_{m-1} \quad \text{ and } \quad q_{m} =  \frac{\chi \, d}{a_0 d + q_{m-1}}  \, q_{m-1} \, ,
\end{align*}
which implies that
\begin{align*}
	r_{m} > r_{m-1} \text{ if } r_{m-1} < d \, \Big(\frac{\chi}{a_0}-1\Big) \quad \text{ and } \quad q_m > q_{m-1} \text{ if } q_{m-1} < d \, (\chi-a_0) \, .
\end{align*}
Due to the choice of $\chi$, these if-conditions are satisfied for $m = 1$. Since $r \mapsto \frac{\chi \, d\, r}{a_0 \, (r + d)}$ and $q \mapsto  \frac{\chi \, d\, q}{a_0 d + q}$ are increasing, we thus can show that $\{r_m\}_{m\in \mathbb{N}_0}$ and $\{q_m\}_{m\in \mathbb{N}_0}$ are strictly increasing sequences with limits 
\begin{align}\label{Ronnie}
	r_m \nearrow d \, \Big(\frac{\chi}{a_0}-1\Big) \quad \text{ and } \quad  q_{m} \nearrow d \, (\chi-a_0) \, .
\end{align}
For further use, we note the following consequence: By means of the properties $\chi < 1$, $a_0 = d /(d+\bar{r})$ and the definition of $\bar{r}$, we see that
\begin{gather}\label{CON4}
	r_m < \frac{m_0 \, d}{d+ (\beta_1-1) \, m_0} \quad \text{ and } \quad 1 + \frac{d}{r_m} - \frac{d}{m_0} > \beta_1 \, , \qquad q_{m} < \frac{d\bar{r}}{d+\bar{r}} \leq s_0 \, .
\end{gather}
We next fix $0 < R_0$ and we let $R_m = R_0 \, (1+2^{-m})/2$ for $m \in \mathbb{N}$. We let $G_m = \Omega \cap B_{R_m}(0)$. We have $G_m \subset\!\!\subset G_{m-1}$. Applying the inequality \eqref{GN} for $G = G_{m} \cap \Omega_i$, $i = 1,2$, we obtain that
\begin{align}\label{KOLFF}
		\|\nabla u\|_{L^{r_m}(G_m)} \leq C_1(m) \, \|D^2 u\|_{L^{a_0r_m}(G_m\setminus S))}^{a_0} \, [u]_{C^{\lambda_0}(G_m)}^{1-a_0} + C_2(m) \, [u]_{C^{\lambda_0}(G_m)} \, ,
\end{align}
Moreover, we have by definition that $a_0 \, r_m = q_{m-1} = \chi \, r_{m-1} \, d/(r_{m-1}+d)$. Since $\chi < 1$, we also have that $q_{m-1} <  r_{m-1} \, d/(r_{m-1}+d)$. Recall that we verified in \eqref{CON4} that $1 + d/r_m - d/m_0 > \beta_1$ independently on $m$. The assumptions of Lemma \ref{SecondDeriv} are therefore valid with $q = q_{m-1}$ and $r = r_{m-1}$, which implies that
\begin{align*}
\|D^2 u\|_{L^{q_{m-1}}(G_m\setminus S))} = & \|D^2 u\|_{L^{a_0r_m}(G_m\setminus S))} \\
\leq & C(m) \, (\|u\|_{W^{1,r_{m-1}}(G_{m-1})} + \|E_{\beta_1} f\|_{L^{s_1}(G_{m-1})} + \|f\|_{L^{q_{m-1}}(G_{m-1})}) \, .
\end{align*}
We invoke \eqref{CON4}$_3$ and H\"older's inequality to further obtain that
\begin{align}\label{KOLFF2}
	\|D^2 u\|_{L^{q_{m-1}}(G_m\setminus S))} \leq C(m) \, (\|u\|_{W^{1,r_{m-1}}(G_{m-1})} + \|E_{\beta_1} f\|_{L^{s_1}(G_{m-1})} + \|f\|_{L^{s_0}(G_{m-1})}) \, .
\end{align}
Combining the inequalities \eqref{KOLFF} and \eqref{KOLFF2}, elementary arguments allow to deduce that
	\begin{align}\label{UN}
	\|\nabla u\|_{L^{r_m}(G_m)} \leq & C(m) \, (\|\nabla u\|_{L^{r_{m-1}}(G_{m-1})} + \|E_{\beta_1} f\|_{L^{s_1}(G_{m-1})} + \|f\|_{L^{s_0}(G_{m-1})} + \|u\|_{C^{\lambda_0}(G_{m-1}) })\nonumber\\
	\leq & C(m) \, (\|\nabla u\|_{L^{r_{m-1}}(G_{m-1})} + \|E_{\beta_1} f\|_{L^{s_1}(G_{m-1})} + \|f\|_{L^{s_0}(G_{m-1})}+ \|u\|_{L^{s_0^*}(G_{m-1}) })	 \, ,
	\end{align}
where we also invoked \eqref{unash} to estimate the H\"older norm. Similarly, we can apply Lemma \ref{SecondDeriv} to show that
\begin{align*}
	\|D^2 u\|_{L^{q_{m}}(G_m\setminus S))} \leq C(m) \, (\|u\|_{W^{1,r_{m}}(G_{m-1})} + \|E_{\beta_1} f\|_{L^{s_1}(G_{m-1})} + \|f\|_{L^{s_0}(G_{m-1})} ) \, .
\end{align*}
The inequality \eqref{GN} implies that
\begin{align*}
	\|\nabla u\|_{L^{r_m}(G_{m-1})} \leq C_1(m) \, \|D^2 u\|_{L^{q_{m-1}}(G_{m-1}\setminus S))}^{a_0} \, [u]_{C^{\lambda_0}(G_{m-1})}^{1-a_0} + C_2(m) \, [u]_{C^{\lambda_0}(G_{m-1})} \, ,
\end{align*}
and combination of both with \eqref{unash} yields
\begin{align}\label{DEUX}
	& \|D^2 u\|_{L^{q_{m}}(G_m\setminus S))}\nonumber\\
	& \leq C(m) \, (\|D^2 u\|_{L^{q_{m-1}}(G_{m-1}\setminus S)}+ \|E_{\beta_1} f\|_{L^{s_1}(G_{m-1})} +\|f\|_{L^{s_0}(G_{m-1})}+ \|u\|_{L^{s_0^*}(G_{m-1}) }) \, .
	\end{align}
Since for $m = 1$ the right-hand sides of \eqref{UN} and \eqref{DEUX} are finite due to \eqref{ugroeg} and Corollary \ref{SecondDerivcoro}, we can conclude.	
\end{proof}
 In the cusp case (A2), the domain $\Omega_1$ below the cusp satisfies the cone condition, so that the inequality \eqref{GN} remains valid for domains of the form $G = \Omega_1\cap \Omega^{\prime}$. However, we were not able to find directly applicable results from the literature for the regularity in $\Omega_2$. Therefore we state the following substitute.
\begin{prop}\label{GNcusp}
	Assume that $\Omega \setminus S = \Omega_1 \cup \Omega_2$, where $S$ is the graph of $\sigma$ satisfying the properties (A1), (A2), (A3). Then, for $0 \leq \lambda_0 < 1$, and $\frac{1-\lambda_0}{2-\lambda_0} \leq a \leq 1$, $r = q/a$ and $q$ is subject to the restriction 
	\begin{align*}
	1 \leq q <\frac{\theta+d-1}{a \, (1-\lambda_0)} \, ,
	\end{align*}
the inequality \eqref{GN} is valid for every $G = \Omega_2 \cap \Omega^{\prime}$ with constants $C_1$ and $C_2$ depending only on $G$ and the parameters.
\end{prop}
\begin{proof}
	In the proof we replace the domain $\Omega_2$ by its intersection with the cylinder $\{(\bar{x},x_d) \, : \, |\bar{x}| < \sigma^{-1}(x_d), \, 0 < x_d < z\}$ where $z >0$ is fixed.  Let $0 < x_d < z$ be fixed and $\rho := \sigma^{-1}(x_d)$. Let further $B^{\prime}_\rho(x_d) = \{(\bar{x},x_d)\, : \, |\bar{x}| \leq \rho\}$ be the $(d-1)-$dimensional ball with radius $\rho$ and centre $(0,x_d)$.
	We at first can choose arbitrary $1 \leq q < + \infty$, $r = q/a$ and $(1-\lambda_0)/(2-\lambda_0) \leq a \leq 1$ if $q \neq d-1$, $(1-\lambda_0)/(2-\lambda_0) \leq a < 1$ if $q = d-1$. On the $(d-1)$-dimensional flat surface $B^{\prime}_{\rho}(x_d)$, the classical Gagliardo-Nirenberg inequality (\cite{nirenberg66b}, Th.\ 1' or Appendix, Lemma \ref{scaleGN}) yields
	\begin{align*}
		\|\nabla_{\bar{x}} u(\cdot,x_d)\|_{L^r(B_\rho^{\prime}(x_d))} \leq & C_1 \, \rho^{(2-\lambda_0)a +\lambda_0-1} \|D^2_{\bar{x}}u(\cdot,x_d)\|_{L^{q}(B_\rho^{\prime}(x_d))}^{a} \, [u(\cdot,x_d)]^{1-a}_{C^{\lambda_0}(B_\rho^{\prime}(x_d))} \\ &+ C_2 \, \rho^{\frac{d-1}{r}+\lambda_0-1} \, [u(\cdot,x_d)]_{C^{\lambda_0}(B_\rho^{\prime}(x_d))} \, ,
	\end{align*}
	where $C_1$ and $C_2$ are two constants depending only on $d$ and $r$. As we assume that $\rho = \sigma^{-1}(x_d)$, integrating twice the differential inequality of assumption (A3) implies that $\rho \geq K_0^{-1}\, x_d^{1/\theta}$. Moreover, with the function $\phi_0$ from \eqref{phinull}, we know that $x_d = \sigma(\rho) \geq \phi_0(R_0) \, \rho$, and thus
	\begin{align*}
		\|\nabla_{\bar{x}} u(\cdot,x_d)\|_{L^r(B_\rho^{\prime}(x_d))} \leq & C_1 \, \Big(\frac{x_d}{\phi_0(R_0)}\Big)^{(2-\lambda_0)a + \lambda_0-1} \,  \|D^2_{\bar{x}}u(\cdot,x_d)\|_{L^{q}(B_\rho^{\prime}(x_d))}^{a} \, [u]^{1-a}_{C^{\lambda_0}(\Omega_2)} \\
		& + C_2 \, \Big(\Big(\frac{x_d}{\phi_0(R_0)}\Big)^{[\frac{d-1}{r}+\lambda_0-1]^+}+ \Big(\frac{x_d^{1/\theta}}{K_0}\Big)^{[\frac{d-1}{r}+\lambda_0-1]^-}\Big) \, [u]_{C^{\lambda_0}(\Omega_2)} \, ,
	\end{align*}
in which $[\cdot]^{-} = \min\{0, \cdot\}$.
	We raise the latter to the power $r$, integrate for $x_d$ over $(0,z)$, and we get
	\begin{align}\label{primehr}
		\|\nabla_{\bar{x}} u\|_{L^r(\Omega_2)}^r \leq  & \tilde{C}_1 \, z^{r((2-\lambda_0)a + \lambda_0-1)}  \|D^2_{\bar{x}}u\|_{L^{q}(\Omega_2)}^{ar} \, [u]^{(1-a)r}_{C^{\lambda_0}(\Omega_2)}\nonumber \\
		& + \tilde{C}_2 \, \Big(z^{1+[d-1-(1-\lambda_0)r]^+}+  \int_0^z x_d^{[d-1-(1-\lambda_0) \, r]^-/\theta} \, dx_d\Big)  \, [u]_{C^{\lambda_0}(\Omega_2)}^r \, .
	\end{align}
	The integral $\int_0^z x_d^{[d-1-(1-\lambda_0) \, r]^-/\theta} \, dx_d$ is finite under the condition 	$r < (\theta+d-1)/(1-\lambda_0)$. Thus, we obtain the estimate \eqref{GN} for the transversal component  $\nabla_{\bar{x}}$ of the gradient.

	In order to obtain a similar estimate for the $x_d-$derivatives, we apply the same Gagliardo-Nirenberg inequality to the function $u(\bar{x}, \cdot)$ over the interval $I(\bar{x}) = ]\, \sigma(|\bar{x}|), \, 2z\, [$, where we fix $0 < |\bar{x}| < \sigma^{-1}(z) =: \rho_1$. Hence
	\begin{align*}
		\|\partial_{x_d} u(\bar{x},\cdot)\|_{L^r(I(\bar{x}))} \leq & C_1 \, |I(\bar{x})|^{(2-\lambda_0)a +\lambda_0-1} \|D^2_{x_d}u(\bar{x},\cdot)\|_{L^{q}(I(\bar{x}))}^{a} \, [u(\bar{x},\cdot)]^{1-a}_{C^{\lambda_0}(I(\bar{x}))} \\ &+ C_2 \, |I(\bar{x})|^{\frac{1}{r}+\lambda_0-1} \, [u(\bar{x},\cdot)]_{C^{\lambda_0}(I(\bar{x}))} \, ,
	\end{align*}
	where $C_1$ and $C_2$ are two constants depending only on $r$. Since $z \leq |I(\bar{x})|  \leq 2z$, these terms do not create any singularity. Thus, raising the latter to the power $r$, we obtain after integration over $B_{\rho_1}^{\prime}(\sigma(\bar{x}))$ that
		\begin{align*}
		\|\partial_{x_d} u\|_{L^r(\Omega_2)} \leq & \tilde{C}_1 \, \|D^2_{x_d}u\|_{L^{q}(\Omega_2)}^{a} \, [u]^{1-a}_{C^{\lambda_0}(\Omega_2))} + \tilde{C}_2 \,  \, [u]_{C^{\lambda_0}(\Omega_2)} \, .\qedhere
	\end{align*}
	\end{proof}	
Then we can state the main result of this section.
\begin{prop}
	Under the assumptions of Lemma \ref{Lippeschitz} we define
	\begin{align*}
		r_2 := d \, \min\Big\{\frac{m_0}{[d+(\beta_1-1)\, m_0]^+}, \, \frac{1+(\theta-1)/d}{[1-\lambda_0]^+} \Big\}\, .
	\end{align*}
	Then $\nabla u \in L^r_{\text{loc}}(\Omega_2; \, \mathbb{R}^d)$ for all $1 \leq r < r_2$ and $D^2 u \in L^q_{\text{loc}}(\Omega_2; \, \mathbb{R}^{d\times d})$ for all $1 \leq q < r_2d/(r_2+d)$. For all $\Omega^{\prime\prime}\subset\!\!\subset \Omega^{\prime}\subset\!\!\subset\Omega$, there is a constant $c = c(r,q,\Omega^{\prime\prime},k_0,k_1)$ such that
	\begin{align*}
		\|D^2 u\|_{L^{q}(\Omega^{\prime\prime}_2)} + 	\|\nabla u \|_{L^{r}(\Omega^{\prime\prime}_2)} \leq c \, (\|E_{\beta_1} f\|_{L^{s_1}(\Omega^{\prime})} +\|f\|_{L^{s_0}(\Omega^{\prime})} + \|u\|_{L^{s_0^*}(\Omega^{\prime})})
	\end{align*} 
	%
	%
\end{prop}
\begin{proof}
We argue as in the proof of Lemma \ref{Lippeschitz}.	
\end{proof}
Combining the argument of Lemma \ref{Lippeschitz} for the regularity in $\Omega_1$ and the latter statement for the regularity in the inward cusp domain $\Omega_2$, we can regard the proof of Theorem \ref{main2} as complete.

\section{ Sources on the interface}\label{sources}

In this section we finally consider the case $Q \neq 0$. We consider the problem \eqref{P1}, \eqref{P2}, but this time with $f = 0$. Then the weak formulation reads $\int_{\Omega} \kappa \, \nabla u \cdot \nabla \phi \, dx + \int_S Q \, \phi \, da = 0$. The surface measure on $S$ is the $(d-1)-$ dimensional Hausdordff--measure $\mathscr{H}^{d-1}$, and we have abbreviated $da = d\mathscr{H}^{d-1}$.

In order to obtain estimates like \eqref{ugroeg}, \eqref{unash} we must at first determine under which conditions the functional 
\begin{align}\label{A0}
	a^0(Q, \, \phi) := \int_S Q \, \phi \, da \quad \text{ for } \quad \phi \in C_c^1(\Omega)  \, , 
\end{align}
extends an element of $[W^{1,p^{\prime}}_0(\Omega)]^*$. At second, in order to obtain the higher regularity, the relation corresponding to \eqref{v1distribution} reads\footnote{Assume here for simplicity that $\kappa$ is constant.}
\begin{align}\label{v1distributionbis}
	\int_{\Omega} \kappa \nabla v_{\epsilon} \cdot \nabla \phi  \, dx = \int_{\Omega} D_{\kappa}(\tau^{\epsilon}) \nabla u \cdot \nabla \phi \, dx + \int_{S} Q\, \tau^{\epsilon} \cdot \nabla \phi \, da \quad \text{ for all } \quad \phi \in C^1_c(\Omega)\, ,
\end{align}
which motivates determining the conditions on a function $P$ defined on $S$, for which
\begin{align}\label{A1}
	a^1_{\delta}(P, \, \phi) := \int_S |x|^{\delta} \, P \, \tau \cdot \nabla \phi \, da  \quad \text{ for } \quad \phi \in C_c^1(\Omega), \, 0 \leq \delta \leq 1  
\end{align}
extends to an element of $[W^{1,p^{\prime}}_0(\Omega)]^*$. In order to obtain these estimates, we begin by introducing a localisation procedure to decompose $S$ into Lipschitz surfaces. The proof of the following statement is given in the appendix. We denote by $\partial B_1^{d-1}(0)$ be the unit sphere in $d-1$ dimensions.
\begin{prop}\label{localise}
 We can find $m \in \mathbb{N}$, a partition of unity $\eta_0, \, \eta_1, \ldots, \eta_m$ on $\partial B_1^{d-1}(0)$ and relatively open hypersurfaces $S_0, \, S_1, \ldots, S_m \subset S$ enjoying the properties
\begin{enumerate}[(a)]
\item for $i = 0, \ldots,m$, $S_i$ is a manifold of class $\mathscr{C}^{0,1}$;
\item for $x \in S$ the function $x \mapsto \eta_i(\bar{x}/|\bar{x}|)$ vanishes on $S \setminus S_i$;
\item $S \subset \bigcup_{i=0}^m \overline{S_i}$.	
\end{enumerate}	
\end{prop}
For $Q\,:  S \rightarrow \mathbb{R}$, we then define
\begin{align}\label{Qloc}
	Q_i(x)  := \eta_i\Big(\frac{\bar{x}}{|\bar{x}|}\Big) \, Q(x) \quad \text{ for } \quad x \in S\,, i= 0,\ldots,m \, .
\end{align}
Using the properties of partitions of unity, we can represent $Q(x) = \sum_{i=0}^{m} Q_i(x)$ for $x \in S$. Note that each $Q_i$ is supported in the surface $S_i$. Let us also preliminarily remark that for $\sigma$ satisfying (A1) and (A2) or (A2')
\begin{align}\label{remarkintegr}
	\int_S |x|^{-\lambda} \, da < + \infty \quad \text{ for all } \quad 0 \leq \lambda < d-1 \, ,
\end{align}
which we prove in the Appendix, Lemma \ref{xhochminuslambda}. We next estimate $a^0(Q, \cdot )$ defined via \eqref{A0} and $a^1_{\delta}(P, \cdot )$ defined via \eqref{A1}. 
\begin{lemma}\label{bilinear}
If $Q \in L^{p^{\rm \#}_S}(S) \supset {\rm tr}_S(W^{1,p^{\rm \#}}(\Omega))$, then $a^0(Q, \cdot )$ defined via \eqref{A0} extends to a element of $[W^{1,p^{\prime}}_0(\Omega)]^*$ and
\begin{align*}
	|a^0(Q, \, \phi)|\leq c \, \|Q\|_{ L^{{p^{\rm \#}_S}}(S) } \, \|\phi\|_{W^{1,p^{\prime}}(\Omega)} \, . 
\end{align*}
Let $\delta >\max\{0,1-\frac{d-1}{p}\}$. If $P\in {\rm tr}_S(W^{1,p}(\Omega))$, then  $a^1_{\delta}(P, \cdot )$ defined via \eqref{A1} extends to a element of $[W^{1,p^{\prime}}_0(\Omega)]^*$ and
	\begin{align*}
		|a^1_{\delta}(P,\, \phi)| \leq ({\rm diam}(\Omega))^{\frac{\delta}{p}} \, (1+c_{\delta})^{1-\frac{1}{p}} \, \|P\|_{{\rm tr}_S(W^{1,p}(\Omega))} \, \|\phi\|_{W^{1,p^{\prime}}(\Omega)} \, .
	\end{align*}
\end{lemma}
\begin{proof}
	For $\phi \in C^1_c(\Omega)$ we have
	\begin{align*}
		\Big|\int_S Q \, \phi \,da\Big| = \Big|\sum_{i=0}^m \int_{S_i} \, Q_i \, \phi \, da\Big| \leq \sum_{i= 0}^m \|Q_i\|_{L^{p^{\rm \#}_S}(S)} \, \| \phi\|_{L^{p^{\rm \#}_S/(p^{\rm \#}_S-1)}(S_i)} \, , 
	\end{align*}	
	where we made use of the triangle and H\"older's inequalities. Since $S_i$ is of class $\mathscr{C}^{0,1}$, the Sobolev space $W^{1,p^{\prime}}(\Omega)$ possesses traces on $S_i$. We let $c_i$ be a bound for the norm of the trace operator as a map from $W^{1,p^{\prime}}(\Omega)$ into $L^{p^{\rm \#}_S/(p^{\rm \#}_S-1)}(S_i)$, and we obtain that
	\begin{align*}
		\Big|\int_S Q \, \phi \,da\Big| \leq \sum_{i= 0}^m c_i \,  \|Q_i\|_{L^{p^{\rm \#}_S}(S_i)} \, \| \phi\|_{W^{1,p^{\prime}}(\Omega)} \, . 
	\end{align*}
	In turn, the Sobolev embedding theorem implies that $W^{1,p^{\rm \#}}$ embedds into $L^{p^{\rm \#}_S}(S_i)$, and $\|Q\|_{L^{p^{\rm \#}_S}(S_i)} \leq c_i \, \|Q\|_{\text{tr}_S(W^{1,p^{\rm \#}})}$.	Hence, $Q \mapsto a^0(Q,\cdot)$ generates a bounded linear operator on $\text{tr}_S(W^{1,p^{\rm \#}})$ with values in $[W^{1,p^{\prime}}_0(\Omega)]^*$ as claimed. \\
	
		For $i = 0, \ldots,m$, and $0 < \delta \leq 1$, and $P$ defined on $S$, we next consider the bilinear form $a^1_{i\delta}(P,\phi) := \int_{S_i} \eta_i\, P \, |x|^{\delta} \, \tau \cdot \nabla \phi \, da$. On the one hand, using that the $\eta_i$s are bounded by one and $\tau$ is a unit vector field, H\"older's inequality implies that
	\begin{align}\label{interpo1}
		|a^1_{i\delta}(P,\, \phi)| \leq R_0^{\delta} \, \|P\|_{L^p(S_i)} \, \|\nabla_S \phi\|_{L^{p^{\prime}}(S_i)} \, .
	\end{align}
	On the other hand, if $P \in W^{1, \, p}(S_i)$, then using integration by parts, together with the facts that $\eta_i = 0$ on $\partial S_i$ and $\tau \cdot \nabla \eta_i = 0$,\footnote{This follows from the fact that $\eta_i = \eta_i(\bar{x}/|\bar{x}|)$.} it follows that
	\begin{align}\label{C0}
		a^1_{i\delta}(P,\, \phi) = - a_{
			i\delta}(\phi,P) - \int_{S_i} \eta_i \, P \, \phi \, \divv_S( |x|^{\delta} \, \tau) \,  da \, .
	\end{align}
	Using the Lemma \ref{tauregu}, we bound $|\divv_S (|x|^{\delta}\,\tau) | \leq c \, |x|^{-1+\delta}$. We next distinguish three cases. If $p < d-1$, then with $p^* = (d-1)p/(d-1-p)$, H\"older's inequality yields
	\begin{align}\label{C1}
		\Big|\int_{S_i} \eta_i \, P \, \phi \, \divv_S (|x|^{\delta} \, \tau) \, da\Big| \leq \|P\|_{L^{p^*}(S_i)} \, \|\phi\|_{L^{p^{\prime}}(S_i)} \, \||x|^{-1+\delta}\|_{L^{d-1}(S_i)} \, .
	\end{align}
For $\delta > 0$, we have $\||x|^{-1+\delta}\|_{L^{d-1}(S_i)} < + \infty$ owing to \eqref{remarkintegr}. At second, if $p > d-1$, then we get
	\begin{align}\label{C2}
		\Big|\int_{S_i} \eta_i \, P \, \phi \, \divv_S (|x|^{\delta} \, \tau) \, da\Big| \leq \|P\|_{L^{\infty}(S_i)} \, \|\phi\|_{L^{p^{\prime}}(S_i)} \, \||x|^{\delta-1}\|_{L^{p}(S_i)}
	\end{align}
If $\delta > 1- (d-1)/p$, then again $\||x|^{\delta-1}\|_{L^{p}(S_i)} < + \infty$ owing to \eqref{remarkintegr}. The third case is $p = d-1$. Then, we fix a number $\chi \in ]1-\delta, \, 1[$, and with $s$ such that $(d-1)/s = 1-(1-\delta)/\chi$ we get  
	\begin{align}\label{C3}
		\Big|\int_{S_i} \eta_i \, P \, \phi \, \divv_S (|x|^{\delta} \, \tau) \, da\Big| \leq \|P\|_{L^{s}(S_i)} \, \|\phi\|_{L^{p^{\prime}}(S_i)} \, \||x|^{\delta-1}\|_{L^{\frac{\chi \, (d-1)}{1-\delta}}(S_i)}
	\end{align}
	Combining \eqref{C1}, \eqref{C2} and \eqref{C3} with the Sobolev embedding theorem for the $(d-1)-$dimensional Lipschitz--graph $S_i$, we get $|\int_{S_i}\eta_i\,  P \, \phi \, \divv_S (|x|^{\delta} \, \tau) \, da\Big| \leq c_{i\delta} \, \|P\|_{W^{1,p}(S_i)} \, \|\phi\|_{L^{p^{\prime}}(S_i)}$, and hence \eqref{C0} shows that
	\begin{align}\label{interpo2}
		|a^1_{i\delta}(P,\phi)| \leq (1+c_{i\delta}) \, \|P\|_{W^{1,p}(S_i)} \, \| \phi\|_{L^{p^{\prime}}(S_i)} \, . 
	\end{align}
We recall \eqref{interpo1} and \eqref{interpo2}, and applying interpolation (see \cite{lionsmagenesIII}) between $L^p(S_i)$ and $W^{1,p}(S_i)$, we get
\begin{align}\label{interpo3}
|a^1_{i\delta}(P,\phi)| \leq R_0^{\frac{\delta}{p}} \, (1+c_{i\delta})^{1-\frac{1}{p}} \, \|P\|_{W^{1-1/p,p}(S_i)} \, \|\phi\|_{W^{1-1/p^{\prime},p^{\prime}}(S_i)} \, .
\end{align}
It remains to observe that $a^1_{\delta} = \sum_{i=0}^m a^1_{i\delta}$ and, since $S_i$ is a Lipschitz graph and $W^{1,q}(\Omega) \hookrightarrow W^{1-1/q}_q(S_i)$, the latter yields
\begin{align*}
&	|a^1_{\delta}(P,\phi)| \leq c \, R_0^{\frac{\delta}{p}} \, (1+c_{\delta})^{1-\frac{1}{p}} \, \|P\|_{{\rm tr}_S(W^{1,p}(\Omega))} \, \|\phi\|_{W^{1,p^{\prime}}(\Omega)} \, . \qedhere
\end{align*}
\end{proof}
{\bf Proof of the results.} In this paragraph we want to prove Theorem \ref{main} for $Q \neq 0$. We at first need to prove for $Q \neq 0$ the result of Proposition \ref{v1prop}. 
\begin{prop}\label{v1proprand}
	With a constant $p_0$ occurring in \eqref{Meyers}, we consider an arbitrary number $r \in ]\max\{p_0^{\prime},d/(d-1)\}, \, p_0^*[$. Assume that $\|E_{\beta_2}Q\|_{\text{tr}_S(W^{1,s_2}(\Omega_1))} < + \infty$ with parameters $0 < \beta_2 \leq 1$ and $\max\{p_0^{\prime},r^{\rm \#}\} < s_2 \leq + \infty$ subject to the restriction
	\begin{align*}
&		\frac{d}{r} - \frac{1}{\min\{p_0,s_2\}} > 0\\
&		\begin{cases}
			\beta_2 \leq \min\Big\{1, \, \frac{d-1}{r}\Big\} & \text{ for } r \leq \min\{p_0,s_2\}\, ,\\
			\beta_2 < \min\Big\{1 + \frac{d}{r} - \frac{d}{\min\{p_0,s_2\}}, \,	\frac{d}{r} - \frac{1}{\min\{p_0,s_2\}}\Big\} & \text{ otherwise} \, .
		\end{cases}
	\end{align*}	
	Suppose that $u \in W^{1,r}_{\text{loc}}(\Omega)$ is a local weak solution to \eqref{P1}, \eqref{P2}, where $f = 0$. Then the function $v_1 := \tau \cdot \nabla u$ belongs to $\bigcap_{1\leq q < \frac{dr}{d+r}} W^{1,q}_{\text{loc}}(\Omega)$. For each $ 1\leq q < r^{\rm \#}$ and $\Omega^{\prime\prime} \subset\!\!\subset \Omega^{\prime} \subset\!\!\subset \Omega$ compactly
	there is $c = c(r,q,s_2,\beta_2,\Omega^{\prime},k_0,k_1)$ such that
	\begin{align*}
		\|\nabla v_1\|_{L^{q}(\Omega^{\prime\prime})} \leq c \, (\|\nabla u\|_{L^{r}(\Omega^{\prime})} + \|E_{\beta_2} Q\|_{\text{tr}_S(W^{1,s_2}(\Omega^{\prime}))}) \, .
	\end{align*} 
\end{prop}
\begin{proof} 
With $g_{\epsilon}(|x|) = f_{\epsilon}(|x|^{\beta})/|x|^{\beta_2+\delta}$ with $\delta>0$, we have
\begin{align*}
	\int_{S} Q \, \tau_{\epsilon}\cdot \nabla \phi \, da= a^1_{\delta}\Big(g_{\epsilon}(|x|)\,  E_{\beta_2}Q, \, \phi\Big) 	\, .
\end{align*}
For $1 <s < +\infty$ and $\delta > \max\{0,1-(d-1)/s\}$, the Lemma \ref{bilinear} implies that
\begin{align*}
	\Big|	\int_{S} Q \, \tau_{\epsilon}\cdot \nabla \phi \, da\Big|
	\leq C(\delta,s) \, \|g_{\epsilon} \, E_{\beta_2}Q \|_{\text{tr}_S(W^{1,s}(\Omega))} \, \|\phi\|_{W^{1,s^{\prime}}(\Omega)} \, .
\end{align*}
If we restrict to $\beta > \beta_2 + \max\{0,1-(d-1)/s\}$ and $\max\{0,1-(d-1)/s\} < \delta < \beta - \beta_2$, we have
\begin{align*}
|g_{\epsilon}(|x|)| \leq |x|^{\beta-\beta_2 - \delta} \leq R_0^{\beta-\beta_2 - \delta}\quad \text{ and } \quad  |\nabla g_{\epsilon}| \leq \beta \, |x|^{\beta-\beta_2 - \delta - 1} \, ,
\end{align*}
and this implies with $\epsilon_0:= \beta-\beta_2 - \delta$ that
\begin{align*}
	\|g\|_{L^{\infty}(\Omega)} \leq  R_0^{\epsilon_0}, \quad \|\nabla g_{\epsilon}\|_{L^d(\Omega)} \leq \beta \, \||x|^{\epsilon_0-1}\|_{L^d(\Omega)} \, .
\end{align*}
Let $1 < s < d$, and consider an extension $P$ of class $W^{1,s}(\Omega)$ for $E_{\beta_2} Q$. We can bound
\begin{align*}
	\|g_{\epsilon} \, P\|_{L^{s}(\Omega_1)} \leq & R_0^{\epsilon_0} \, \|P \|_{L^{s}(\Omega_1)} \, ,\\
	\|\nabla (g_{\epsilon} \, P) \|_{L^{s}(\Omega_1)} \leq &  R_0^{\epsilon_0} \,\|\nabla  P \|_{L^{s}(\Omega_1)} + \|\nabla g_{\epsilon}\|_{L^{d}(\Omega_1)} \, \|P\|_{L^{\frac{ds}{d-s}}(\Omega)}\\
	\leq & C\, (R_0^{\epsilon_0} + \beta \, \||x|^{-1+\epsilon_0}\|_{L^{d}(\Omega_1)}) \, \|P\|_{W^{1,s}(\Omega_1)} \, .
\end{align*}
Under the conditions $E_{\beta_2}Q \in {\rm tr} (W^{1,s}(\Omega))$ it therefore follows for arbitrary $\beta > \beta_2 + \max\{0,1-(d-1)/s\}$ and $1 \leq s < d$ 
 \begin{align*}
 	\Big|	\int_{S} Q \, \tau_{\epsilon}\cdot \nabla \phi \, da\Big| \leq 
 	C(\beta,\beta_2) \, \|E_{\beta_2}Q \|_{\text{tr}_S(W^{1,s}(\Omega))} \, \|\phi\|_{W^{1,s^{\prime}}(\Omega_1)} \, .
 \end{align*}
Invoking the property \eqref{Meyers}, for $\Omega^{\prime\prime} \subset\!\!\subset \Omega^{\prime} \subset\!\!\subset \Omega$, for $p_0^{\prime} < s < \min\{p_0,d\}$ and for 
\begin{align}\label{con2prime}
\beta > \beta_2+ \max\{0,1-(d-1)/s\} \quad \text{ and } \quad s \leq s_2 \, ,
\end{align}
we show that $v_{\epsilon}$ subject to \eqref{v1distributionbis} satisfies
\begin{align*}
	\|\nabla v_{\epsilon}\|_{L^{s}(\Omega^{\prime\prime})} \leq & c \, (\|D_{\kappa}(\tau^{\epsilon}) \nabla u\|_{L^{s}} + \|a^1_{\beta}(g_{\epsilon} \, E_{\beta_2}Q,\cdot)\|_{[W^{1,s^{\prime}}_0]^*} +\||\tau^{\epsilon}|\,  \nabla u\|_{L^{s}}  )\\
	  \leq & c \, (\|D_{\kappa}(\tau^{\epsilon}) \nabla u\|_{L^{s}} + \|E_{\beta_2}Q\|_{\text{tr}_S(W^{1,s_2}(\Omega))} +\||\tau^{\epsilon}|\,  \nabla u\|_{L^{s}} ) \, .
\end{align*}
Note that the condition \eqref{con2prime} replaces the condition \eqref{CON2} in the proof of Proposition \ref{v1prop}, hence the proof that we can find suitable parameters is more complicated in the present case.  In fact \eqref{con2prime} introduces the restriction $s < \frac{d-1}{1-\beta + \beta_2}$, which is sharper than $\beta \geq \beta_2$.

In order to be able choose $s$ arbitrarily close to the desired threshold $s \nearrow \frac{dr}{d + (1-\beta) \, r}$ with suitable $\beta$ as in the proof of Proposition \ref{v1prop}, we require the inequality 
\begin{align*}
	\frac{dr}{d + (1-\beta) \, r} \leq \frac{d-1}{1-\beta + \beta_2} \, \quad \text{ equivalent with } \quad \beta \geq 1- d \, \Big(\frac{d-1}{r}- \beta_2\Big) \, ,
\end{align*}
which is a new condition for $\beta$.
Hence, we let
\begin{align}\label{tildebeta2}
	\tilde{\beta}_2 = \max\Big\{\beta_2, \,  1-d \, \Big(\frac{d-1}{r}-\, \beta_2\Big)\Big\} \, .
\end{align}
Now the algebraic calculation of the proof of Proposition \ref{v1prop} can be repeated one to one with $(\tilde{\beta}_2,s_2)$ replacing $(\beta_1,s_1)$. With $m_0 = \min\{p_0,s_2\}$, we choose $\beta$ and $s$ according to (cf.\ \eqref{CON3}, \eqref{CON3bis})
\begin{align*}
	&	\begin{cases}
		\tilde{\beta}_2 \leq \beta \leq 1 & \text{ if }  m_0 \geq r \text{ and }  \tilde{\beta}_2 > 1 - d \, (\frac{1}{p_0^{\prime}}-\frac{1}{r}) \, ,\\
		1 - d \, (\frac{1}{p_0^{\prime}}-\frac{1}{r}) < \beta \leq 1 & \text{ if } m_0 \geq r  \text{ and }  \tilde{\beta}_2 \leq  1 - d \, (\frac{1}{p_0^{\prime}}-\frac{1}{r})\, ,\\
		\max\{\tilde{\beta}_2,\,  1 - d \,(\frac{1}{p_0^{\prime}}-\frac{1}{r})\} < \beta <  1 + \frac{d}{r} - \frac{d}{m_0} & \text{ if } m_0 < r \, ,\\
	\end{cases}\\
 & p_0^{\prime} < s < \frac{dr}{d + (1-\beta) \, r}  \, 	.\qedhere
\end{align*}
\end{proof}
For the case $f = 0$, we thus have the following statement (cp.\ Corollary \ref{v1coro}).
\begin{coro}\label{v1coroQ}
	Let $\max\{p_0^{\prime},d/(d-1)\} < p < p_0 $ with the constant $p_0$ occurring in \eqref{Meyers}. We assume that $Q \in L^{p^{\rm\#}_S}(S)$ with $E_{\beta_0}Q \in B^{p,p}_{1-1/p}(S)$, $\beta_0 = \min\{1,(d-1)/p\}$ and that $u$ is a local weak solution to \eqref{P1}, \eqref{P2} satisfying $u \in L^{p}_{\text{loc}}(\Omega)$. Then the function $v_1 := \tau \cdot \nabla u$ belongs to $\bigcap_{1\leq q < p^{\rm \#}} W^{1,q}_{\text{loc}}(\Omega)$. For each $ 1\leq q < p^{\rm \#}$ and $\Omega^{\prime\prime} \subset\!\!\subset \Omega^{\prime} \subset\!\!\subset  \Omega$ compactly, there is $c = c(p,q,\Omega^{\prime\prime},k_0,k_1)$ such that
	\begin{align*}
		\|\nabla v_1\|_{L^{q}(\Omega^{\prime\prime})} \leq c \, (\|Q\|_{L^{p^{\rm\#}_S}(S \cap \Omega^{\prime})} + \|E_{\beta_0}Q\|_{B^{p,p}_{1-1/p}(S \cap \Omega^{\prime})}+\|u\|_{L^p(\Omega^{\prime})}) \, .
	\end{align*} 
\end{coro}
This Corollary is proved applying the preceding proposition with $r = p$, $s_2 = p$, and choosing $\beta_2 = \min\{1,(d-1)/p\}$. In order to prove the Theorem \ref{main} with $Q \neq 0$, we observe that (for $f = 0$ and $Q \neq 0$ sufficiently regular) the statement of Lemma \ref{tangentes} must read
\begin{align*}
	& \int_{\Gamma_R} \kappa \, \diff u \cdot \diff \zeta \, dS +  \int_{\Gamma_R \cap S} (v_1 \, [\kappa]_S\tau \cdot \nu+Q) \, \zeta\, ds  = \int_{\Gamma_R} \Phi \, \zeta \, dS \quad \text{ for all } \quad \zeta \in W^{1,2}(\Gamma_R) \, .
\end{align*}
Thus, repeating for $Q$ the estimates done for the term $v_1 \, [\kappa]_S\tau \cdot \nu$, we obtain for all $1 < s < +\infty$
\begin{align*}
	\|D^2 u\|_{L^s(\Gamma_R \setminus S)} \leq C \, (&\|f\|_{L^s(\Gamma_{R})} + \|\nabla v_1\|_{L^s(\Gamma_{R})} + \|\nabla Q\|_{L^s(\Gamma_{R})} \\
	& + \frac{1}{R} \, (\|\nabla u\|_{L^s(\Gamma_{R})} + \|Q\|_{L^s(\Gamma_{R})} )+ \frac{1}{R^2} \, \|u\|_{L^s(\Gamma_{R})}) \, ,
\end{align*}
with a constant independent on $u$ and for all $0< R \leq R_0$. Now we can finish the proof of Theorem \ref{main} exactly as in the case $Q = 0$.

\appendix

\section{Auxiliary estimates}

Let $n \geq 1$ denote the space dimension. For $R > 0$, we denote by $B_R$ the ball of radius $R$ centred at zero in $\mathbb{R}^n$. We consider fields $\psi$ defined on $\partial B_1$ and $\Psi$ defined in $B_R$, where we assume $R \geq R_1 > 1$. Further we let $A \in C^{0,1}(B_R \setminus \partial B_1; \, \mathbb{R}^{n\times n})$ be a matrix value map, such that $A(x)$ is elliptic for all $x \in B_R\setminus \partial B_1$, with uniform ellipticity constant $a_0>0$.
\begin{lemma}\label{TransmiBasic}
	Let $1 < s < + \infty$. Let $\Psi \in L^{s}(B_R)$ and $\psi \in W^{1-1/s}_s(\partial B_1)$. Suppose that $w \in W^{1,2}_0(B_R)$ satisfies
\begin{align*}
	\int_{B_R} A(x) \nabla w \cdot\nabla \zeta \, dx = \int_{B_R} \Psi \, \zeta \, dx + \int_{\partial B_1} \psi \, \zeta \, da \quad \forall \zeta \in C^1(\overline{B_R}) \, .
\end{align*}
Then, there is $c = c(s,R_1,n)$ such that $\|\nabla w\|_{W^{1,s}(B_R)} \leq c (\|\Psi\|_{L^s(B_R)} + \|\psi\|_{W^{1-1/s}_s(\partial B_1)})$.
\end{lemma}
\begin{proof}
For the $W^{2,s}$ regularity, we refer to the Theorem 6.5.2 of \cite{MR3524106}, and we obtain $\|w\|_{W^{2,s}(B_R)} \leq c_R (\|\Psi\|_{L^s(B_R)} + \|\psi\|_{W^{1-1/s}_s(\partial B_1)})$. We next want to show that the constant can be chosen independent of $R$.

Since the testfunctions are not constrained and $A(x)$ is regular, we get $A(x)\nabla w \cdot \nu = 0$ on  $\partial B_R$, and since $w = 0$ on $\partial B_R$ by assumption, it follows that $\nabla w = 0$ on $\partial B_R$ in the sense of traces. Hence the extension by zero of $w$ outside of $B_R$ belongs to $W^{2,s}(\mathbb{R}^n)$ and satisfies
\begin{align}\label{nondiv}
-A(x) \, : \, D^2 w - (\nabla \cdot A(x)) \cdot \nabla w = \Psi(x) \, \chi_{B_R} \text{ in } \mathbb{R}^n\, , \quad [\![A(x)\nabla w \cdot \nu ]\!] = \psi \text{ on } \partial B_1 \, , 
\end{align}
where the operator on the left-hand side is independent on $R$. Invoking Theorem 6.5.2 of \cite{MR3524106} and the main result of the paper \cite{MR0298209}, we obtain separately uniform estimates in $W^{2,s}(B_{R_1})$ (for the bounded domain case with interface) and in $W^{2,s}(\mathbb{R}^n\setminus B_{1+(R_1-1)/2})$ (for the exterior problem with continuous coefficients). Alternatively to applying the result of \cite{MR0298209}, we can replace the operator in \eqref{nondiv} by $-A(x) \, : \, D^2 w - (\nabla \cdot A(x)) \cdot \nabla w + w$ and solve the equation with a modified right-hand side $\Psi + w$. The $L^s$ norm of $\Psi + w$ over the whole space can be handled by means of the Lemma \ref{transmibasicorderlow} here below. Then, the Theorem 6.1.10 in  \cite{MR3524106}, shows that the modified operator associated with the exterior problem is $\mathcal{R}$-sectorial and possesses a bounded inverse from $L^s(\mathbb{R}^n)$ into $W^{2,s}$ for all $1 < s < + \infty$. 
%
%
%
%
\end{proof}
\begin{lemma}\label{transmibasicorderlow}
	Under the assumptions of Lemma \ref{TransmiBasic}, there exists for all $1 < s < +\infty$ a constant $c = c(n,s,a_0)$ independent on $R$ such that
	\begin{align*}
		\|w\|_{L^{s} (\mathbb{R}^n)} \leq c \, (\|\Psi\|_{L^s(B_R)} + \|\psi\|_{W_s^{1-1/s}(\partial B_1)})\, . 
	\end{align*}
\end{lemma}
\begin{proof}
We perform the proof for $n >2$ and let the case $n = 2$ as an exercise.
Assume that $w$ is a sufficiently regular solution to 
\begin{align*}
	\int_{\mathbb{R}^n} A(x) \nabla w \cdot\nabla \zeta \, dx = \int_{\mathbb{R}^n} f_R  \, \zeta \, dx + \int_{\partial B_1} f_S \, \zeta \, da \quad \forall \zeta \in C^1(\mathbb{R}^n) \, ,
\end{align*}
in which $(f,f_S)$ are generic data in $B_R$ and on $\partial B_1$ and $f_R  := f \, \chi_{B_R}$. If $(f,f_S)$ belongs to $L^{2^{\rm \#}}(B_R) \times  L^{2^{\rm \#}_S}(\partial B_1)$, then inserting $w$ as testfunction yields
\begin{align*}
	a_0 \, \int_{\mathbb{R}^n} |\nabla w|^2 \, dx \leq \|f_R\|_{L^{2^{\rm \#}}} \, \|w\|_{L^{2^*}} + \|f_S\|_{L^{2^{\rm \#}_S}} \, \|w\|_{L^{2^*_S}}\, .
\end{align*}
With the trace operator we further have $\|w\|_{L^{2^*_S}} \leq c_0 \, \|w\|_{W^{1,2}(B_1)}$. Since $\|w\|_{L^2(B_1)} \leq |B_1|^{\frac{2}{n}} \, \|w\|_{L^{2^*(B_1)}}$, we obtain the inequality
\begin{align*}
\|w\|_{W^{1,2}(B_1)} \leq c_1 \, \|w\|_{L^{2^*}} + \|\nabla w\|_{L^2} \leq (c_1 \, c_2+1) \, \|\nabla w\|_{L^2} \, ,
\end{align*}
in which $c_2 = c_2(n)$ is the Sobolev constant on $\mathbb{R}^n$. Thus
\begin{align*}
\|\nabla w\|_{L^2} \leq \frac{C}{a_0} \, (\|f_R\|_{L^{2^{\rm \#}}} +\|f_S\|_{L^{2^{\rm \#}_S}}) \, ,
\end{align*}
which by the Sobolev inequality implies that
\begin{align*}
  \|w\|_{L^{2^*}} \leq \frac{C}{a_0} \, (\|f_R\|_{L^{2^{\rm \#}}} +\|f_S\|_{L^{2^{\rm \#}_S}}) \quad \text{ and } \quad \|w\|_{L^{2^*_S}(\partial B_1)} \leq \frac{C}{a_0} \, (\|f_R\|_{L^{2^{\rm \#}}} +\|f_S\|_{L^{2^{\rm \#}_S}}) \, .
\end{align*}
Thus, the solution operator, denoted by $T$, is bounded from $L^{2^{\rm \#}}(B_R) \times  L^{2^{\rm \#}_S}(\partial B_1)$ into $L^{2^*}(B_R)\times  L^{2^{*}_S}(\partial B_1)$, independently of $R$. 

For $(f, \, f_S) \in L^1(B_R) \times L^1(\partial B_1)$ we use the testfunction ${\rm sign}(w) \, (1 - 1/(1+|w|)^\delta)$ where $0< \delta < 1$ is a parameter, to obtain that
\begin{align*}
\delta \,  a_0\, 	\int_{\mathbb{R}^n} \frac{|\nabla w|^2}{(1+|w|)^{1+\delta}} \, dx \leq \|f_R\|_{L^1} + \|f^S\|_{L^1} \, .
\end{align*}
Now with the function $v := {\rm sign}(w) \, ((1+|w|)^{\alpha}-1)$, $\alpha := (1-\delta) /2$, we obtain that
\begin{align*}
	\frac{4\delta}{(1-\delta)^2} \, a_0\,	\int_{\mathbb{R}^n} |\nabla v|^2 dx \leq \|f_R\|_{L^1} + \|f^S\|_{L^1} \, ,
\end{align*}
which also yields
\begin{align*}
	\|v\|_{L^{2^*}}^2 \leq \frac{C_\delta}{a_0} \, (\|f_R\|_{L^{1}} +\|f_S\|_{L^{1}}) 
	\quad \text{ and } \quad \|v\|_{L^{2^*_S}}^2 \leq \frac{C_\delta}{a_0} \, (\|f_R\|_{L^{1}} +\|f_S\|_{L^{1}})  \, .
\end{align*}
We can use the inequality $2^{1-\alpha} \, |v| \geq |w|^{\alpha}$ to obtain that
 \begin{align*}
 	\|w\|_{L^{\alpha \, 2^*}}^{2\alpha} \leq \frac{C_\delta}{a_0} \, (\|f_R\|_{L^{1}} +\|f_S\|_{L^{1}})  	\quad \text{ and } \quad \|v\|_{L^{\alpha 2^*_S}(\partial B_1)}^{2\alpha} \leq \frac{C_\delta}{a_0} \, (\|f_R\|_{L^{1}} +\|f_S\|_{L^{1}})  \, .
 \end{align*}
Note that $\alpha$ attains every value in $(0,1/2)$ for some $\delta$ in $(0,1)$. Thus, we see that for every $q \in [1,2^*/2[$ and $q_S \in [1, \, 2^*_S/2[$, the operator $T$ is bounded from $L^1(B_R) \times L^1(\partial B_1)$ into $L^q(B_R) \times L^{q_S}(\partial B_1)$ independently of $R$.

With the Riesz-Thorin interpolation theorem we obtain after some computations that for all $1 \leq p \leq 2^{\rm\#}$ and all $1 \leq p_S \leq 2^{\rm \#}_S$ the operator $T$ is bounded independently of $R$ from $L^{p}(B_R) \times L^{p_S}(\partial B_1)$ into $L^{r}(B_R) \times L^{r_S}(\partial B_1)$, where the parameter $r$ is subject to the restriction
\begin{align*}
	1 \geq  \frac{1}{r} > \frac{1}{p} - \frac{2}{n} \quad \text{ and } 1 \geq \frac{1}{r} > \Big(1-\frac{1}{n}\Big) \frac{1}{p_S} - \frac{1}{n}  \, . 
\end{align*} 
We can calculate that for $1 < s \leq 2^{\rm \#}$, the inequalities are valid with $r = s$, $p = s$ and $p_S = s_S^* = s(n-1)/(n-s)$. Hence
\begin{align*}
	\|w\|_{L^{s}} \leq \frac{C_\delta}{a_0} \, (\|f_R\|_{L^{s}} +\|f_S\|_{W^{1-1/s}_s(\partial B_1)}) \, .
\end{align*}
Using the fact that $T$ is selfadjoint in $L^2(B_R) \times L^2(\partial B_1)$, we can also obtain the boundedness independently on $R$ for all $s \in (1, +\infty)$ via a duality argument. We spare in this place the technical details.
\end{proof}

%
%
%
%

The next Lemma states the scaling behaviour for the inequality \eqref{GN}.
\begin{lemma}\label{scaleGN}
	Let $q$, $\lambda$ and $a$ be constants subject to the restrictions \eqref{RESTRI1}. We denote by $C_1$ and $C_2$ two constants such that \eqref{GN} is valid with $G = B_1 \subset \mathbb{R}^n$. Then, for $\rho >0$, \eqref{GN} holds true on $B_{\rho}$ in the form
	\begin{align}\label{GNrho}
		\|\nabla u\|_{L^{q/a}} \leq C_1 \, \rho^{(2-\lambda)a + \lambda -1}\,  \|D^2 u\|_{L^{q}}^a \, [u]_{C^{\lambda}}^{1-a} + C_2 \, \rho^{\lambda - 1 + a \, \frac{n}{q}}\, [u]_{C^{\lambda}} \, ,
	\end{align}
	where all norms are taken over the domain $B_{\rho}$.
\end{lemma}
\begin{proof}
	For $u$ defined on $B_{\rho}$ we define a new function $\tilde{u}(z) := u (\rho \, z)$ for $z \in B_1$. Applying \eqref{GN} to $\tilde{u}$, the claim follows by means of the transformation formula. In particular, notice that $\|D^{\alpha}_z \tilde{u}\|_{L^q(B_1)} = \rho^{|\alpha| -n/q} \, \|D^{\alpha} u\|_{L^q(B_{\rho})}$ and $[\tilde{u}]_{C^{\lambda}(B_1)} = \rho^{\lambda} \, [u]_{C^{\lambda}(B_{\rho})}$.
\end{proof}

%

{\bf Proof of Proposition \ref{localise}.} We first consider the case $d \geq 3$. For $i = 0, \ldots, 16$, we denote by $J_i$ the interval $](i-1)\pi/8,\, \, (i+1)\pi/8[$. Using the spherical coordinates 
\begin{align*}
	\xi = \begin{pmatrix}
		\cos(\phi_{1})\\
		\sin(\phi_{1}) \, \cos(\phi_{2})\\
		\vdots \\
		\sin(\phi_{1})\cdots \sin(\phi_{d-3})\cos(\phi_{d-2})\\
		\sin(\phi_{1})\cdots \sin(\phi_{d-3})\sin(\phi_{d-2})
	\end{pmatrix}\, ,
\end{align*}
a point $\xi$ on the $(d-1)-$dimensional unit sphere $\partial B_1^{d-1}(0) =: S^{d-1}_1$ is represented by means of $d-2$ angles $\phi_1(\xi) ,\ldots,\phi_{d-2}(\xi) \in P^{d-2} := ]0,\pi[ \times \ldots \times ]0,\pi[ \times ]0,2\pi]$. We can cover the parallelepiped $P^{d-2}$ by means of
\begin{align*}
	P^{d-2} \subset \bigcup_{k_{1} = 0}^{8} \ldots \bigcup_{k_{d-3} = 0}^{8} \bigcup_{k_{d-2}=0}^{16} J_{k_1} \times \ldots \times J_{k_{d-3}} \times J_{k_{d-2}} \, .
\end{align*}
With $m := (16)^{d-2}/2^{d-3}$ and appropriately ordering the multi-indices $(k_1,\ldots,k_{d-2})$ we define $U_i = \{\xi \in S^{d-1}_1 \, :\, \phi(\xi) \in J_{k_1(i)} \times \ldots \times J_{k_{d-3}(i)} \times J_{k_{d-2}(i)}\}$ for $i = 0, \ldots,m$. In this way, we cover $	S^{d-1}_1 = \bigcup_{i=0}^{m} U_i$. Then, we let $\eta_0, \ldots, \eta_{m}$ be a corresponding partition of unity. This means that $\eta_i \in C_c^{\infty}(U_i)$ is nonnegative and $\sum_{i=0}^{m} \eta_i\equiv 1$ on $S_1^{d-1}$. 

Now, we cover the surface $S$ by means of the pieces 
\begin{align*}
	S_i = \{x \in S \, : \, \phi_1(\bar{x}/|\bar{x}|),\ldots, \phi_{d-2}(\bar{x}/|\bar{x}|) \in U_i   \} \, .
\end{align*}
We claim that $S_i$ is a Lipschitz--graph. To see the latter, denoting by $\xi$ a reference point on $U_i$, we have $|\phi(\bar{x}/|\bar{x}|) - \phi(\xi)| \leq \pi/8$ on $S_i$. For $x \in S_i$, it is readily shown that
\begin{align*}
	|\bar{x}| = \bar{x} \cdot \Big(\frac{\bar{x}}{|\bar{x}| } - \xi\Big) + \bar{x} \cdot \xi \leq \frac{\pi}{8} \,  |\bar{x}|+ |\bar{x} \cdot \xi | \, ,  
\end{align*}
which implies that $|\bar{x} \cdot \xi| \geq 	|\bar{x}| \, (1-\pi/8) $ and, since $x_d = \sigma(|\bar{x}|)$, that
\begin{align*}
	|\bar{x} - (\bar{x}\cdot\xi) \, \xi| \leq \sqrt{\frac{\pi}{16-\pi}} \, \sigma^{-1}(x_d)\quad \text{ for all } \quad x \in S_i \, .
	%
\end{align*}
We find that
\begin{align*}
	S_i = \{x \in S \, : \, \bar{x}/|\bar{x}| \in U_i,\, \bar{x} \cdot \xi = \sqrt{(\sigma^{-1}(x_d))^2 - |\bar{x}- (\bar{x}\cdot\xi) \, \xi|^2 }  \} \, .
\end{align*}
This implies the inclusion $S_i \subset \text{graph}(f)$ in the coordinate system with the $d^{\rm th}$-axis given by $(\xi, \, 0)$, where $f$ is the function $(t_1, \, \ldots, \, t_{d-2}, \, z) \mapsto \sqrt{(\sigma^{-1}(z))^2-|t|^2}$ defined in the open set $G := \{(t_1, \ldots,t_{d-2},\, z) \, : \, 0 < z < \sigma(R_0), \, |t| <  \sqrt{\pi/(16-\pi)} \, \sigma^{-1}(z)\}$. The gradient of $f$ is easily shown to be bounded by a constant times $1 + \max_{0 \leq z \leq \sigma(R_0)} |(\sigma^{-1})^{\prime}(z)|$.\\

{\bf Case $d=2$:} 
This case is particular because we cannot localise the surface using partitioning into sectors. We define $m = 1$, $S_0 = S \cap \{x \in \mathbb{R}^2 \, : \, x_1 > 0\}$ and $S_1 := S \cap \{x \in \mathbb{R}^2 \, : \, x_1 < 0\}$. For $Q: \, S \rightarrow \mathbb{R}$, we define $Q_i = Q \, \eta_i({\rm sign}(x_1))$ with $\eta_0(t) = (1+t)/2$ and $\eta_1(t) = (1-t)/2$.

\begin{lemma}\label{xhochminuslambda}
Assume that $\sigma$ satisfies the conditions (A1) and (A2). Then $\int_S |x|^{-\lambda} \, da < + \infty$ for all $0 \leq \lambda < d-1$.
\end{lemma}	
\begin{proof}
We first estimate $\int_S |x|^{-\lambda} \, da \leq 2^{\lambda/2} \,  \int_S 1/(|\bar{x}|^{\lambda}+\sigma(|\bar{x}|)^{\lambda}) \, da$. Thus
\begin{align}\label{flat0}
	\int_S \frac{1}{|x|^{\lambda}} \, da \leq & 2^{\lambda/2} \,  \int_{B_{R_0}^{d-1}(0)} \frac{1+\sigma^{\prime}(|\bar{x}|)}{|\bar{x}|^{\lambda}+\sigma(|\bar{x}|)^{\lambda}} \, d\bar{x} \nonumber \\
	=  & 2^{\lambda/2} \, \omega_1 \, \int_0^{R_0} r^{d-2} \, \frac{1+\sigma^{\prime}(r)}{r^{\lambda}+\sigma(r)^{\lambda}} \, dr \, ,
\end{align}
where $\omega_1$ is the measure of the unit sphere in $d-1$ dimensions. If $d = 2$, then for $\lambda < d-1 = 1$, we readily obtain that
\begin{align*}
	 \int_0^{R_0} \frac{1+\sigma^{\prime}(r)}{r^{\lambda}+\sigma(r)^{\lambda}} \, dr \leq \frac{1}{1-\lambda} \, (R_0^{1-\lambda}+ \sigma(R_0)^{1-\lambda}) \, 
\end{align*}
and the claim follows from \eqref{flat0}. Thus, we discuss the case $d > 2$ and assume without loss of generality that $d-1 > \lambda >1$. Then the first contribution on the right-hand of \eqref{flat0} obviously satisfies 
\begin{align*}
\int_0^{R_0} r^{d-2} \, \frac{1}{r^{\lambda}+\sigma(r)^{\lambda}} \, dr \leq \frac{R_0^{d-1-\lambda}}{d-1-\lambda} \, .
\end{align*}
In order to estimate the second contribution, we let $\psi(r) = \sigma(r)/r$, and we note that
\begin{align*}
	\int_0^{R_0} r^{d-2} \, \frac{\sigma^{\prime}(r)}{r^{\lambda}+\sigma(r)^{\lambda}} \, dr = \int_0^{R_0} r^{d-1-\lambda} \, \frac{\psi^{\prime}(r)}{1+\psi(r)^{\lambda}} \, dr  + \int_0^{R_0} r^{d-2-\lambda} \, \frac{\psi(r)}{1+\psi(r)^{\lambda}} \, dr \, .
\end{align*}
As $\lambda > 1$, the second contribution is again bounded by $R_0^{d-1-\lambda}/(d-1-\lambda)$. For the first contribution, let us define for $t > 0$ the function $\phi(t) = \int_0^{t} 1/(1+s^{\lambda}) \, ds$, which is globally bounded. We have
\begin{align*}
&	\int_0^{R_0} r^{d-1-\lambda} \, \frac{\psi^{\prime}(r)}{1+\psi(r)^{\lambda}} \, dr = R_0^{d-1-\lambda} \, \phi(R_0) - (d-1-\lambda)\, \int_0^{R_0} r^{d-2-\lambda} \, \phi(r) \, dr \, . \qedhere
	\end{align*}
\end{proof}


\end{document}